\documentclass{article}

\usepackage[T1]{fontenc}
\usepackage[latin1]{inputenc}
\usepackage{epsfig}
\usepackage{amsmath,amssymb}
\usepackage[a4paper]{geometry}
\usepackage{graphicx}
\usepackage{caption}
\usepackage{subcaption}
\def\XXint#1#2#3{{\setbox0=\hbox{$#1{#2#3}{\int}$}
    \vcenter{\hbox{$#2#3$}}\kern-.5\wd0}}

\def\longrightharpoonup{\DOTSB\relbar\joinrel\rightharpoonup}

\def\RR{\mathbb R}
\def\ZZ{\mathbb Z}
\def\NN{\mathbb N}

\def\11{\mathbf 1}
\begin{document}
\numberwithin{equation}{section}
\newtheorem{theoreme}{Theorem}[section]
\newtheorem{proposition}[theoreme]{Proposition}
\newtheorem{remarque}[theoreme]{Remark}
\newtheorem{lemme}[theoreme]{Lemma}
\newtheorem{corollaire}[theoreme]{Corollary}
\newtheorem{definition}[theoreme]{Definition}

\title{Homogenization of heat diffusion in a cracked medium}
\author{Xavier Blanc$^1$, Benjamin-Edouard Peigney$^2$ \\
\\
{\footnotesize $^1$ Laboratoire Jacques-Louis Lions,} \\
{\footnotesize Université Paris-Diderot,} \\
{\footnotesize 75205 PARIS Cedex 13, FRANCE.} \\
{\footnotesize {\tt blanc@ann.jussieu.fr}} \\
\\
{\footnotesize $^2$ CEA, DAM, DIF,}\\
{\footnotesize 91297 ARPAJON Cedex, FRANCE.}\\
{\footnotesize {\tt Benjamin.Peigney@cea.fr}}}

\maketitle

\begin{abstract}
We develop in this note a homogenization method to tackle the problem of a diffusion process
through a cracked medium. We show that the cracked surface of the domain
induces a source term in the homogenized equation. We assume that the
cracks are orthogonal to the surface of the material, where an incoming
heat flux is applied. The cracks are supposed to be of depth 1, of small
width, and periodically arranged.
\end{abstract}

\section{Motivation and setting of the problem}
\label{sec:pb}

\begin{figure}[h!]
  \centering
  \includegraphics[width=0.7\linewidth]{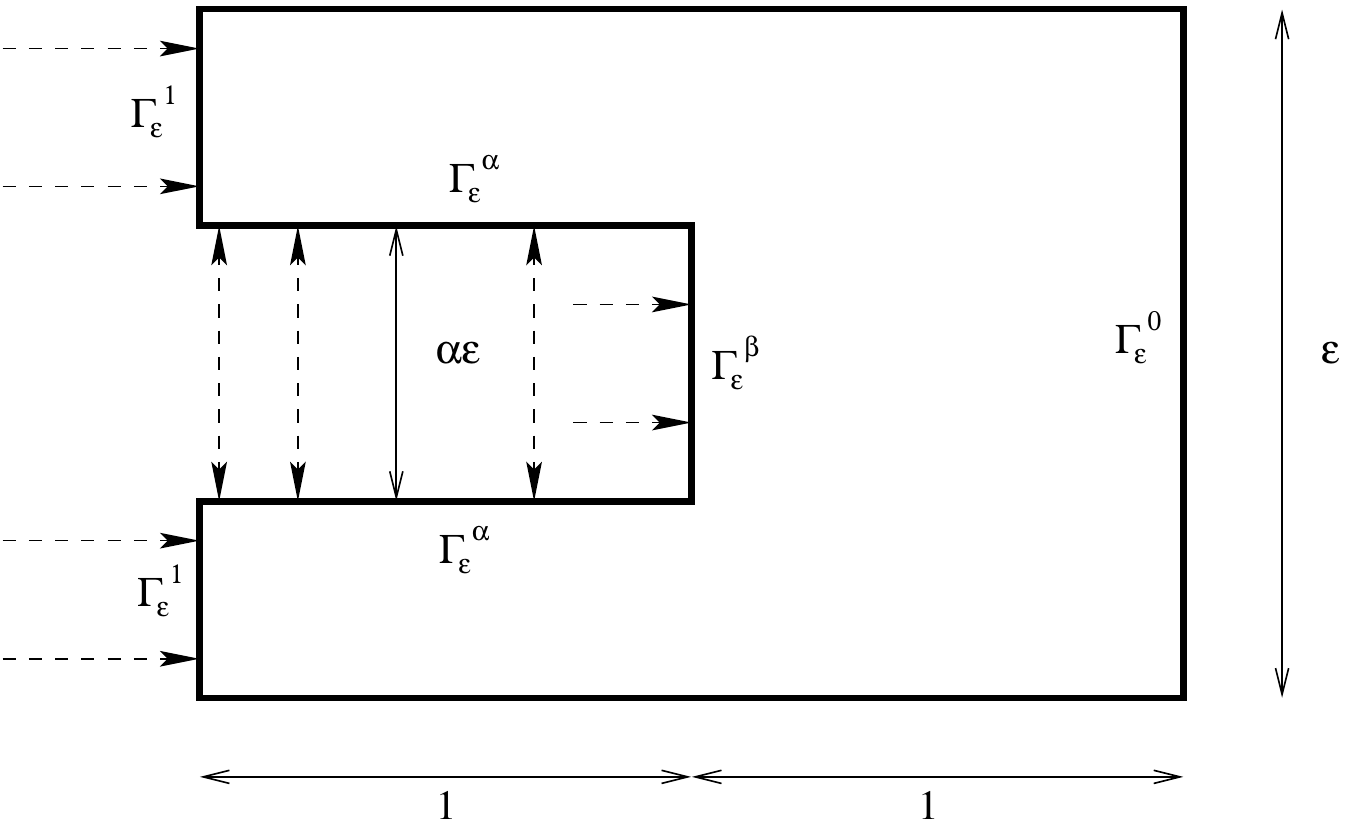}
  \caption{The cracked domain $\Omega_\varepsilon$. \label{fig:dessin}}
\end{figure}

We consider the propagation of radiation through a cracked medium, made of an optically thick material. The propagation is initiated by an incoming given energy flux imposed on the left boundary of the crack (see figure~\ref{fig:dessin} where the flux is represented in dashed lines).

Physically, the exchange surface between the optically thick medium and the source may be greatly modified by the fractures. This may have a significant impact on the energy balance of the considered system. In many situations, we cannot model the surface of the cracked medium directly, which is often too intricate to be described. Besides, the shape of the fractures may have a stochastic feature and it may involve many spatial scales. Full numerical simulations of such multi-scaled media become hence infeasible.

That is why we have been looking for an \textbf{average} approach, to capture the effects of cracks in a homogenized medium. The model presented here is simple enough to be coupled to standard FEM codes. The physical idea behind the model developed in this paper, called "MOSAIC" (Method Of Sinks Averaging Inhomogeneous behavior of Cracked media), is to treat the flux enhancement induced by the crack as a volume \textbf{source term} in the homogenized energy equation. We will show that this can be justified rigorously by homogenization theory.

Since the medium is assumed to be optically thick, the propagation of radiation follows a diffusion process \cite{mih99}. For the sake of simplicity, we shall assume a linear behavior law, which means that the energy flux $\mathcal{F}$ is proportional to the energy gradient $\nabla u$.\footnote{But the results given here could be extended to the non-linear (power law) case.}

Besides, we consider the diffusion process on a short time scale, so that the hydrodynamic effects, which are not supposed to play a significant role, are disregarded. 


The linear diffusion problem can thus be modeled by:

\begin{equation}
  \label{eq:simple}
  \begin{cases}
    \partial_t u -\Delta u = 0 & \text{in } \Omega_\varepsilon, \\
    \partial_n u = 0 & \text{on } \Gamma_\varepsilon^0  \\
    \partial_n u = 1 & \text{on } \Gamma_\varepsilon^1, \\
    \partial_n u = \frac{\alpha-\beta} 2 \varepsilon & \text{on }
    \Gamma_\varepsilon^\alpha,\\
    \partial_n u = \frac{\beta}\alpha &\text{on } \Gamma_\varepsilon^\beta.
  \end{cases}
\end{equation}
We also need an initial condition:
\begin{equation}
  \label{eq:ci_simple}
  u(x,y,t=0) = u^0(x,y),
\end{equation}
so that problem (\ref{eq:simple})-(\ref{eq:ci_simple}) is well-posed.
For the problem at hand, $u(x,y,t)$ is the energy density but it could represent any field following a diffusion process.
 
We impose $u(x,y,t)$ to be periodic of period $\varepsilon$ with respect
to the $y$ direction. The domain $\Omega_\varepsilon$, as well as the
boundaries
$\Gamma_\varepsilon^0,\Gamma_\varepsilon^1,\Gamma_\varepsilon^\alpha,
\Gamma_\varepsilon^\beta$ are defined on figure~\ref{fig:dessin}. The
right-most boundary $\Gamma_\varepsilon^0$ is supposed to coincide with
the set $\{x=1\}$, while $\Gamma_\varepsilon^\beta$ is a subset of
$\{x=0\}$, and the left-most boundary $\Gamma_\varepsilon^1$ is a subset
of $\{x=-1\}.$ The
parameter $\varepsilon > 0$ is supposed to be small and will tend to
$0$, whereas $\alpha\in [0,1)$ is a fixed parameter related to the width of the
crack. The parameter $\beta\in[0,\alpha)$ measures the portion of the flux which,
coming through the segment $\{x=-1, -\alpha\varepsilon/2<y<\alpha\varepsilon/2\}$, reaches the
bottom $\Gamma^\beta_\varepsilon$ of the crack. The remaining part of
the incoming flux is distributed on the horizontal part of the boundary,
namely $\Gamma^\alpha_\varepsilon$. The parameter $\beta$ is supposed to
be fixed. The boundary conditions in (\ref{eq:simple}) are defined in
such a way that the total incoming flux is exactly equal to $1$, which
is the value we impose on the left boundary in the case $\alpha=0$ (no crack).


A space-time dependence of the flux applied on the boundaries
$\Gamma_\varepsilon^\alpha$ may be introduced but it does not affect the
homogenization process that we describe here. 


\section{Changing the scale}
\label{sec:scale}

To carry out an asymptotic expansion of the solution $u=u_\varepsilon$ of (\ref{eq:simple}) in powers of $\varepsilon$, we "scale" the variable $y$, in the spirit of \cite{blp}. Actually, 2 scales describe the model: the variable $y$ is the macroscopic one, whereas $\displaystyle\frac y \varepsilon$ represents the "microscopic geometry". Thus, we define:

$$u(x,y,t) = v\left(x, \frac y \varepsilon,t\right),$$
so that $v$ is periodic of period $1$.

We notice that:

$$\partial_xu(x,y,t) = \partial_xv \left(x,\frac y \varepsilon,t\right),
\quad \partial^2_x u(x,y,t) = \partial^2_x v\left(x,\frac y \varepsilon,t\right),$$
$$ \partial_y   u(x,y,t) = \frac 1 \varepsilon\partial_y v \left(x,\frac y
    \varepsilon,t\right),\quad \partial^2_y(x,y,t) = \frac 1
  {\varepsilon^2}\partial^2_y v\left(x,\frac y \varepsilon,t\right).$$
And $v$ is solution of (\ref{eq:simple_v}):

\begin{equation}
  \label{eq:simple_v}
  \begin{cases}
    \displaystyle-\partial^2_x v -\frac 1 {\varepsilon^2}
    \partial^2_y v  + \partial_t v = 0 & \text{ in } \Omega_1, \\
    \partial_n v = 0 & \text{on } \Gamma_1^0, \\
    \partial_n v = 1 & \text{on } \Gamma_1^1, \\
    \displaystyle\frac 1 {\varepsilon^2}\partial_n v = \frac{\alpha-\beta} 2  & \text{on }
    \Gamma_1^\alpha, \\
    \partial_n v = \frac\beta\alpha &  \text{on } \Gamma_1^\beta
  \end{cases}
\end{equation}

\section{Asymptotic expansion}
\label{sec:homog}

\begin{figure}[h!]
  \centering
  \includegraphics[width=0.7\linewidth]{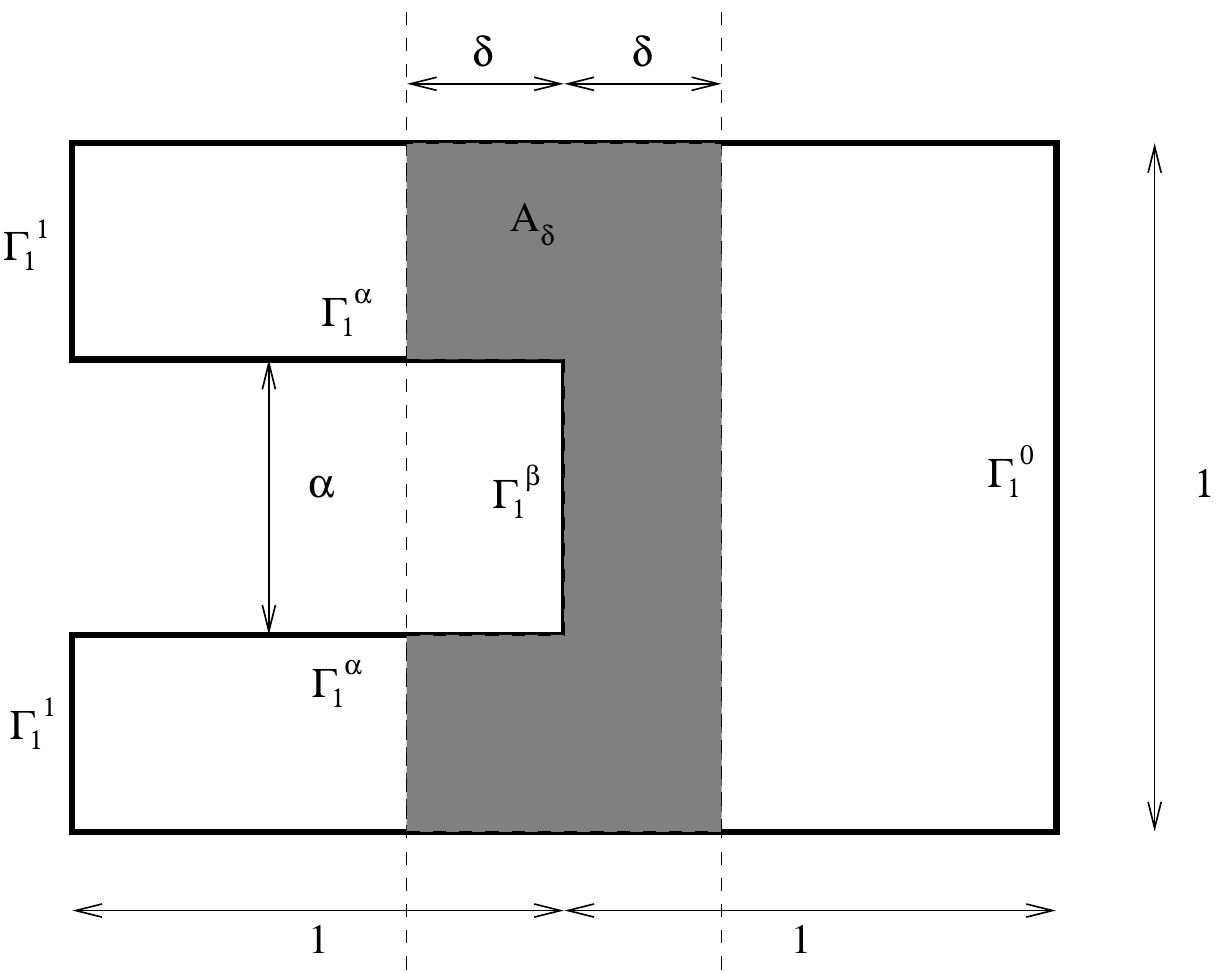}
  \caption{The domain $A_\delta$ on which we integrate.}
  \label{fig:A_delta}
\end{figure}

Firstly, we notice that in the system (\ref{eq:simple_v}) the domain $\Omega_1$ does not depend of $\varepsilon$ anymore. We have to study an equation depending on $\varepsilon$  in a fixed domain. Besides, the parameter $\varepsilon$ appears in the equation (\ref{eq:simple_v}) only as $\varepsilon^2$, which means that $\varepsilon^2$ is a good parameter for an asymptotic expansion. 
Thus, it seems natural to look for $v$ as follows:
\begin{equation}
  \label{eq:asymptotique}
  v(x,y,t) = v_0(x,y,t) + \varepsilon^2 v_1(x,y,t) + \varepsilon^4 v_2(x,y,t) +
  \dots
\end{equation}
Hence, we insert the ansatz (\ref{eq:asymptotique}) into the system
(\ref{eq:simple_v}) and identify the different powers of $\varepsilon^2$.
We obtain:

\begin{itemize}
\item \underline{At the order $\varepsilon^{-2}$ :}
$$\partial^2_y v_0 = 0,$$
which means that $v_0$ can be written
$$v_0(x,y,t) = f(x,t) + y g(x,t),$$
where $f$ and $g$ are 2 functions independent of the variable $y$. The condition of periodicity in $y$ verified by $v$ implies that $v_0(x,1/2,t) = v_0(x,-1/2,t)$, hence: $g=0$.
Thus, $v_0$ does not depend on $y$ :
\begin{equation}
  \label{eq:v_0}
  v_0(x,y,t) = v_0(x,t).
\end{equation} 
Besides, the boundary conditions on $v_0$ are:
\begin{equation}
\label{eq:boundary_v0}
\partial_n v_0 = 0 \text{ on } \Gamma_1^0, \quad \partial_n v_0 = 1
\text{ on } \Gamma_1^1, \quad \partial_n v_0 = 0 \text{ on }
\Gamma_1^\alpha, \quad \partial_n v_0 = \frac\beta\alpha \text{ on } \Gamma_1^\beta.
\end{equation}

We check that the boundary condition on $\Gamma_1^\alpha$ is consistent
with (\ref{eq:v_0}), because it is equivalent to the statement:
$\displaystyle\partial_y v_0 = 0$ on $\Gamma_1^\alpha,$ which is automatically verified since $v_0$ does not depend on $y$.

\item \underline{At the order $\varepsilon^0$ :}

\begin{equation}
\label{eq:dev_1}
-\partial^2_x v_0 - \partial^2_y v_1 + \partial_t v_0 = 0.
\end{equation}

The boundary conditions on $v_1$ give
$$\partial_n v_1 = 0 \text{ on } \Gamma_1^0, \quad \partial_n v_1 = 0
\text{ on } \Gamma_1^1, \quad \partial_n v_1 = \frac{\alpha-\beta} 2 \text{ on }
\Gamma_1^\alpha, \partial_n v_1 = 0 \text{ on } \Gamma_1^\beta.$$

We now integrate (\ref{eq:dev_1}) with respect to $y$ : we get
\begin{equation}
\label{eq:dev_2}
-(1-\alpha) \partial^2_x v_0 + (1-\alpha) \partial_t
v_0 = \int_{-\frac12}^{-\frac\alpha 2} \partial^2_y v_1 dy
+ \int^{\frac12}_{\frac\alpha 2} \partial^2_y v_1 dy,
\end{equation}
if $x<0$, and
$$- \partial^2_x v_0 + \partial_t
v_0 = \int_{-\frac12}^\frac12 \partial^2_y v_1
dy,$$
if $x>0.$ (recall that we assume that the axis $x=0$ contains the border $\Gamma_1^{\beta}$.)
In the first case $x<0$, we use the boundary value for $\partial_n v_1$,which gives
$$\partial_y v_1\left(x,-\frac \alpha 2,t\right) = 
\frac{\alpha-\beta}2, \quad \partial_y v_1\left(x,\frac
  \alpha 2,t\right) = -\frac{\alpha-\beta}2.$$
Besides, the periodicity in $y$ implies that boundary values in $y=1/2$ are exactly compensated by those in $y=-1/2$. Thus, the right hand side of equation (\ref{eq:dev_2}) is $\frac {\alpha-\beta}{2} + \frac{\alpha-\beta}{2} = \alpha-\beta$, and we obtain
$$
\begin{cases}
  \displaystyle -(1-\alpha) \partial^2_x v_0 + (1-\alpha) \partial_t
v_0 = \alpha-\beta & \text{if } x<0, \\
\displaystyle - \partial^2_x v_0 + \partial_t
v_0 = 0 & \text{if } x>0.
\end{cases}
$$
\end{itemize}
Hence, $v_0$ satisfies an equation in $\{x<0\}$, and an equation in
$\{x>0\}$. In order to define it properly, we need boundary
conditions. For $\Gamma_1^0$ and $\Gamma_1^1,$ we have
(\ref{eq:boundary_v0}). We are now going to derive boundary conditions
on $\{x=0\}$. Since we have assumed that $\partial_n u = \beta/\alpha$ on
$\Gamma_1^{\beta}$, and since $v_0$ does not depend on $y$, it may seem
natural to impose that $\partial_x v_0(x=0) = -\beta/\alpha.$ However, as
we will see below, the flux $\partial_x v_0$ is not continuous across
the interface $\{x=0\}$. Therefore, we need to take into account a
boundary layer at this interface. For this purpose, we go back to
(\ref{eq:simple_v}), in which we zoom at $x=0$, that is, we define
$$v(x,y,t) = w\left(\frac x \varepsilon, y,t\right).$$
Inserting this into (\ref{eq:simple_v}), we see that, 
\begin{equation}\label{eq:equation_w}
-\partial_x^2 w - \partial_y^2 w +\varepsilon^2\partial_t w= 0.
\end{equation}
We integrate this equation over the domain $A_\delta = [-\delta,\delta]\times
[-1/2,1/2]\cap \Omega_1$ (see figure~\ref{fig:A_delta}), integrate by
parts, and find
\begin{multline*}
\int_{\{x=\delta\}} \partial_x w - \int_{\{x=-\delta,
  \alpha/2<|y|<1/2\}} \partial_x w - \int_{\{x=0,|y|<\alpha/2\}} \partial_x
  w \\- \int_{\{y=\alpha/2,-\delta<x<0\}}\partial_y w +
  \int_{\{y=-\alpha/2,-\delta<x<0\}}\partial_y w =
  \varepsilon^2\int_{A_\delta} \partial_t w 
\end{multline*}
The right-hand side vanishes as $\varepsilon\to 0$. In the left-hand
side, the fourth and fifth terms are bounded by $\delta$, hence, taking
$\varepsilon\to 0$, then $\delta\to 0$, we infer
$$\int_{\{x=0^+\}} \partial_x w =
\int_{\{x=0^-,\alpha/2<|y|<1/2\}}\partial_x w - \beta.$$
Therefore, we have, going back to $v_0$, and
using the fact that it does not depend on $y$,
\begin{equation}
\label{eq:transmission_flux}
(1-\alpha) \partial_x v_0(0^-, y,t) = \partial_x v_0(0^+,y,t) +\beta.
\end{equation}
We can also compute the relation between $v_0(0^+,y,t)$ and
$v_0(0^-,y,t)$ by multiplying (\ref{eq:equation_w}) by $x$, and
integrating it again over $A_\delta$. We then have
\begin{multline*}
\delta\int_{\{x=\delta\}} \partial_x w +\delta \int_{\{x=-\delta,
  \alpha/2<|y|<1/2\}} \partial_x w 
   - \int_{\{y=\alpha/2,-\delta<x<0\}}x\partial_y w +
  \int_{\{y=-\alpha/2,-\delta<x<0\}}x\partial_y w \\
- \int_{\{x=\delta\}}
  w + \int_{\{x=-\delta, \alpha/2<|y|<1/2\}}w +
  \int_{\{x=0,|y|<\alpha/2\}} w
=
  \varepsilon^2\int_{A_\delta} \partial_t w 
\end{multline*}
Here again, the right-hand side vanishes as $\varepsilon\to 0,$ while
the first line vanishes as $\delta\to 0$. Hence, we find that
$$\int_{\{x=0^+\}} w = \int_{\{x=0^-, \alpha/2<|y|<1/2\}} w + \int_{\{x=0,
  |y|<\alpha/2\}} w.$$
Going back to $v_0$, this implies that $v_0$ is continuous accross the
interface.

We can thus write the system of equation satisfied by $v_0$:
\begin{equation}
\label{eq:dev_3.1}
\begin{cases}
- \partial^2_x v_0 + \partial_t
v_0 = \frac{\alpha-\beta}{1-\alpha} & \text{in } \Omega_1\cap\{x<0\}, \\
\partial_n v_0 = 1 & \text{on } \Gamma_1^1, \\
\partial_n v_0 = 0 & \text{on } \Gamma_1^\alpha, \\
\end{cases}
\end{equation} 
and
\begin{equation}
\label{eq:dev_3.2}
\begin{cases}
- \partial^2_x v_0 + \partial_t
v_0 = 0 & \text{in } \Omega_1\cap\{x>0\}, \\
\partial_n v_0 = 0 & \text{on } \Gamma_1^0, \\
\end{cases}
\end{equation} 
This system is not well-posed, since boundary conditions are missing at
the interface $\{x=0\}$. We thus impose the transmission conditions
(\ref{eq:transmission_flux}), together with the fact that $v_0$ should
be continuous accross the interface:
\begin{equation}
  \label{eq:transmission}
  v_0(x=0^-) = v_0(x=0^+), \quad (1-\alpha) \partial_x v_0(x=0^-)
  = \partial_x v_0(x=0^+) +\beta.
\end{equation}



\section{Homogenized equation}

We have formally shown that

$$u(x,y,t) \approx v_0\left(x,\frac y \varepsilon,t\right),$$
with $v_0$ solution of (\ref{eq:dev_3.1})-(\ref{eq:dev_3.2}). Moreover, if we extend $u$ by
$0$ outside $\Omega_\varepsilon$, we use the fact that, for any function
$f$ which is $1$-periodic with respect to $y$, we have
\begin{equation}\label{eq:lim_faible}
f\left(\frac y \varepsilon\right) \mathop{\longrightharpoonup}^*
\int_{-1/2}^{1/2} f(y)dy,
\end{equation}
in $L^\infty$. Hence, $u$ converges to the average of $v_0$ with respect
to $y$, that is, $(1-\alpha) v_0(x,t)$ if $x<0$, and
$v_0(x,t)$ if $x>0$.
In other words, the limit
equation on $u$ is thus (\ref{eq:dev_3.1}) multiplied by $(1-\alpha)$,
and (\ref{eq:dev_3.2}). Note that the boundary conditions are treated
exactly as the equation, using (\ref{eq:lim_faible}). Hence, the system
satisfied by $u$ reads:

\begin{equation}
\label{eq:hmodel.1}
\begin{cases}
  -\Delta u + \partial_t u =  \alpha-\beta &
  \text{in } \{ -1<x<0\}, \\
  \partial_n u = 1- \alpha & \text{on } \{x=-1\}, \\
  u \text{ is }1-\text{periodic in } y.
\end{cases}
\end{equation}
\begin{equation}
\label{eq:hmodel.2}
\begin{cases}
  -\Delta u + \partial_t u =  0 &
  \text{in } \{ 0<x<1\}, \\
  \partial_n u = 0 & \text{on } \{x=1\}, \\
  u \text{ is }1-\text{periodic in } y.
\end{cases}
\end{equation}
We also have the corresponding transmission condition inherited from
(\ref{eq:transmission}):
\begin{equation}
  \label{eq:htransmission}
  u(x=0^-) = (1-\alpha)u(x=0^+), \quad \partial_x u(x=0^-)
  = \partial_x u(x=0^+) + \beta.
\end{equation}
\begin{remarque}\label{rk:hmodel}
  In the coupling condition (\ref{eq:htransmission}), we have a jump for $u$
  as well as for the flux $\partial_x u$. Therefore, it is {\em a
    priori} not possible to recast
  (\ref{eq:hmodel.1})-(\ref{eq:hmodel.2})-(\ref{eq:htransmission}) into
  a single boundary problem in the domain $\{-1< x < 1\}$ (however, see
  Section~\ref{formulation-faible} below for a formal formulation of the problem). Nevertheless, when
  $\alpha$ is small and $\beta=0$, (\ref{eq:htransmission}) almost
  amounts to impose continuity of $u$ and its derivative accross the
  interface $\{x=0\}$. Therefore, in such a case, the approximate
  problem
\begin{equation}
\label{eq:hmodel}
\begin{cases}
  -\Delta u + \partial_t u =  \alpha\11_{\{x<0\}} &
  \text{in } \{ -1<x<1\}, \\
  \partial_n u = 1- \alpha & \text{on } \{x=-1\}, \\
  \partial_n u = 0 & \text{on } \{x=1\}, \\
  u \text{ is }1-\text{periodic in } y,
\end{cases}
\end{equation}
should give a solution which is close to the solution to
(\ref{eq:hmodel.1})-(\ref{eq:hmodel.2})-(\ref{eq:htransmission}).

\end{remarque}

\begin{remarque}\label{rk:hmodel2}
A spatial dependence on the flux imposed on the boundary $\Gamma_\varepsilon^\alpha$ may be introduced. In this case, and setting $\beta = 0$ for the sake of simplicity, the formulation of the problem
written on the cracked domain $\Omega_\varepsilon$ reads:

\begin{equation}
  \label{eq:simple2}
  \begin{cases}
    \partial_t u -\Delta u = 0 & \text{in } \Omega_\varepsilon, \\
    \partial_n u = 0 & \text{on } \Gamma_\varepsilon^0  \\
    \partial_n u = 1 & \text{on } \Gamma_\varepsilon^1, \\
    \partial_n u = f_\alpha(x) \varepsilon & \text{on }
    \Gamma_\varepsilon^\alpha,\\
    \partial_n u = 0 &\text{on } \Gamma_\varepsilon^\beta.
  \end{cases}
\end{equation}

The function $f_\alpha(x)$ is such that:

$$
\int_{x=-1}^{x=0}f_\alpha(x)dx = \frac{\alpha}{2},
$$ 
ensuring the conservation of energy as $\varepsilon$ varies.

Furthermore, if $f_\alpha(x)$ is such that $f_\alpha(0)=0$, then the singularity at $x=0$ vanishes and the
formulation of the homogenized model is exactly given by:

\begin{equation}
\label{eq:hmodel2}
\begin{cases}
  -\Delta u + \partial_t u =  f_\alpha(x)\11_{\{x<0\}} &
  \text{in } \{ -1<x<1\}, \\
  \partial_n u = 1- \alpha & \text{on } \{x=-1\}, \\
  \partial_n u = 0 & \text{on } \{x=1\}, \\
  u \text{ is }1-\text{periodic in } y.
\end{cases}
\end{equation}

In such a case, the solution $u$ and its gradient $\partial_n u$ are continuous across the interface $x=0$.
\end{remarque}




\section{Analysis of the homogenized problem}
\label{sec:analyse_homog}

In this section, we prove that
problem~(\ref{eq:hmodel.1})-(\ref{eq:hmodel.2})-(\ref{eq:htransmission})
is well posed. Let us first set the notation:
\begin{equation}
  \label{eq:def_omega}
  \Omega =  \left\{(x,y)\in \RR^2, \quad
    -\frac12 <y<\frac12, \ -1<x<1\right\}.
\end{equation}
\begin{equation}
\label{eq:def_omega_+}
\Omega^+ = \Omega \cap \{x>0\} = \left\{(x,y)\in \RR^2, \quad
    -\frac12 <y<\frac12, \ 0<x<1\right\}.
\end{equation}
\begin{equation}
\label{eq:def_omega_-}
\Omega^- = \Omega \cap \{x<0\} = \left\{(x,y)\in \RR^2, \quad
    -\frac12 <y<\frac12, \ -1<x<0\right\}.
\end{equation}
\begin{equation}
\label{eq:def_gamma_-}
\Gamma^1 =  \left\{(x,y)\in \RR^2, \quad
    -\frac12 <y<\frac12, \ x=-1\right\}.
\end{equation}
\begin{equation}
\label{eq:def_gamma_+}
\Gamma^0 =  \left\{(x,y)\in \RR^2, \quad
    -\frac12 <y<\frac12, \ x=1\right\}.
\end{equation}
\begin{equation}
\label{eq:def_gamma_0}
\Gamma^\beta =  \left\{(x,y)\in \RR^2, \quad
    -\frac12 <y<\frac12, \ x=0\right\}.
\end{equation}
\begin{lemme}\label{lm:existence_unicite_homog}
  Assume that $\alpha\in (0,1)$. Then for any $T>0$, problem
  (\ref{eq:hmodel.1})-(\ref{eq:hmodel.2})-(\ref{eq:htransmission}) has a
  unique solution $(u_-,u_+) \in X,$ where
$$X = C\left([0,T],H^1(\Omega^-)\right)\cap
  C^1\left( [0,T], L^2(\Omega^-)\right)
  \times C\left([0,T],H^1(\Omega^+) \right)\cap C^1\left([0,T], L^2(\Omega^+)\right)$$
\end{lemme}
\noindent{\bf Proof:} First, let us point out that, for any $F\in
L^2(\Gamma^\beta)$, the problem 
\begin{equation}
  \label{eq:homog_+}
  \begin{cases}
  -\Delta u + \partial_t u =  0 &
  \text{in } \Omega^+, \\
  \partial_n u = 0 & \text{on } \Gamma^0, \\
  \partial_n u = F & \text{on } \Gamma^\beta, \\
  u \text{ is }1-\text{periodic in } y,\\
  u(t=0) = u_0 & \text{ in } \Omega^+.
\end{cases}
\end{equation}
admits a unique solution $u^+ \in C\left([0,T],H^1(\Omega^+)\right) \cap C^1\left([0,T],L^2(\Omega^+)
\right)$, for any $T>0$. This is easily proved using standard tools
of the analysis of elliptic PDEs. See for instance \cite{evans} for the details.

Next, consider the problem (here, $g\in L^2(\Gamma^\beta)$:
\begin{equation}
  \label{eq:homog_-}
  \begin{cases}
  -\Delta u + \partial_t u =  \alpha - \beta &
  \text{in } \Omega^-, \\
  \partial_n u = 1-\alpha & \text{on } \Gamma^1, \\
  u = g & \text{on } \Gamma^\beta, \\
  u \text{ is }1-\text{periodic in } y, \\
  u(t=0) = u_0 & \text{ in } \Omega^-.
\end{cases}
\end{equation}
Here again, standard theory of elliptic equations allows to prove that
(\ref{eq:homog_+}) has a unique solution $u^-
\in C\left([0,T],H^1(\Omega^-)\right) \cap C^1\left([0,T],L^2(\Omega^-)
\right)$, for any $T>0$ (see \cite{evans}). We now study the following
fixed-point approach: consider an initial guess $F^0 \in L^2(\Gamma^\beta)$,
to which we associate the solution $u^{+,0}$ of
(\ref{eq:homog_+}). Then, define $g^0$ as the trace of $(1-\alpha)u^{+,0}$ on
$\Gamma^\beta$, and solve (\ref{eq:homog_-}) with data $g = g^0$: this
defines $u^{-,0},$ and a new flux $F^1 = \beta-\partial_n u^{-,0}$ in
$L^2(\Gamma^\beta)$. Repeating this procedure, we build a sequence
$(F^n)_{n\in\NN}$ in $L^2(\Gamma^\beta)$, together with the corresponding
solutions $u^{\pm,n}$ and the data $g^n$. In order to prove that this
sequence converges, we are going to prove that the application
$\Phi : L^2([0,T]\times \Gamma^\beta) \to L^2([0,T]\times\Gamma^\beta)$, which associates $F^{n+1}$ to
$F^n$ is a contraction mapping. For this purspose, we define
$$G^n = F^{n+1}-F^n, \ v^{n,\pm} = u^{n+1,\pm} - u^{n,\pm}, \ h^n =
g^{n+1} - g^n.$$
It is clear that $v^{n,+}$ satisfies (\ref{eq:homog_+}) with
$F = G^n$ and $u_0 = 0$. Similarly, $v^{n,-}$ is
the solution to (\ref{eq:homog_-}) with $\alpha-\beta = 0,$ $1-\alpha =
0$, $u_0 = 0$, and $g = h^n$. As a consequence, using the uniqueness for
problem~(\ref{eq:homog_-}), we infer
$ v^{n,-}(x,y,t) = (1-\alpha)v^{n,+}(-x,y,t),$ whence
\begin{equation}
\label{eq:G_n_plus_1}
G^{n+1} = (1-\alpha)G^n.
\end{equation}
This immediately implies that the sequence $(F^n)_{n\in\NN}$ converges
in $L^2\left([0,T]\times \Gamma^\beta\right).$ Denoting by $F$ its limit, it
is clear that the corresponding solution $u^{\pm}$ is a solution to
(\ref{eq:hmodel.1})-(\ref{eq:hmodel.2})-(\ref{eq:htransmission}). 

This proves the existence of a solution. Next, considering the
uniqueness, we assume that we have two solutions $u_1$ and $u_2$ such
that $u_1 \neq u_2$, we necessarily have $\partial_x {u_1}_{|\Gamma^\beta} \neq \partial_x
{u_2}_{|\Gamma^\beta}$, according to the uniqueness of the solution to
(\ref{eq:homog_+}). Hence, defining $F_1 = \partial_x {u_1}_{|\Gamma^\beta}$
and $F_2 = \partial_x {u_2}_{|\Gamma^\beta}$, these functions are fixed points
of the above application $\Phi$. Hence, the above argument implies that
$F_1 = F_2$, which proves uniqueness. \hfill$\Box$

\begin{remarque}
  The above proof is in fact useful for numerical purposes. Indeed, it
  proves that this fixed-point approach always converges. Hence, it may
  be used to compute the solution $u$ to
  (\ref{eq:hmodel.1})-(\ref{eq:hmodel.2})-(\ref{eq:htransmission}). However,
  it should be noted that (\ref{eq:G_n_plus_1}) is an {\em
    equality}. Hence, if $\alpha$ is close to $1$, the convergence will
  be very slow.
\end{remarque}

\section{Proof of convergence}
\label{sec:preuve}

In this section, we give a rigourous proof of the fact that the solution
$u_\varepsilon$ to (\ref{eq:simple}) converges to the solution $u$ to
(\ref{eq:hmodel.1})-(\ref{eq:hmodel.2}) as $\varepsilon\to 0$. 

As a preliminary remark, let us point out that, using standard results
of PDE analysis, one easily proves that (\ref{eq:simple}) has a unique
solution in $C\left([0,T],H^1(\Omega_\varepsilon)\right) \cap C^1\left([0,T],L^2(
\Omega_\varepsilon)\right)$, for any $T>0$ and any $\varepsilon>0$. See for instance \cite{evans} for the details.

\begin{proposition}
\label{pr:cv}
  Let $u_\varepsilon$ be the unique solution to (\ref{eq:simple}). We
  extend it by $0$ outside $\Omega_\varepsilon$, and assume that the initial data
  $u_\varepsilon(t=0)$ is such that 
  \begin{equation}\label{eq:cv_0}
    u_\varepsilon(t=0) \longrightharpoonup u_0 \text{ in } L^2(\Omega).
  \end{equation}
Then, for any $T>0$, we have
  \begin{equation}
    \label{eq:convergence}
    u_\varepsilon \mathop{\longrightharpoonup}_{\varepsilon\to 0} u
    \text{ in } L^2(\Omega\times [0,T]),
  \end{equation}
where $u$ is the unique solution to (\ref{eq:hmodel.1})-(\ref{eq:hmodel.2})-(\ref{eq:htransmission}).
\end{proposition}

\begin{remarque}
  \label{rq;zero}
  We have set $u_\varepsilon=0$ outside $\Omega_\varepsilon$. This
  strategy is physically relevant, since
  $u_\varepsilon$ is a temperature, and the heat transfer is only
  modelled inside $\Omega_\varepsilon$: one may think of the outside of
  $\Omega_\varepsilon$ as the vacuum, or at least a domain which is
  transparent to radiation.
\end{remarque}

\begin{remarque}
  In (\ref{eq:convergence}), we have only a weak convergence. The reason
  for this is the fact that we have extended $u_\varepsilon$ by $0$
  outside $\Omega_\varepsilon$, as it is explained in
  Remark~\ref{rq;zero}, whereas it is positive in
  $\Omega_\varepsilon$, due to the imposed incoming flux. Hence,
  $u_\varepsilon$ qualitatively behaves like a function which is equal
  to $1$ in $\Omega_\varepsilon$, and $0$ outside. In the domain
  $\Omega_1$, this function converges weakly to its average, but does
  not converge strongly in $L^2$.
\end{remarque}

\begin{remarque}
  The convergence (\ref{eq:convergence}) is only local in time ($T$
  cannot be infinite). This is due to the fact that we impose a constant
  incoming flux. Therefore, integrating (\ref{eq:simple}) over
  $\Omega_\varepsilon$, and using an integration by parts, we have
$$\frac d {dt} \int_{\Omega_\varepsilon} u = 1,$$
hence $u_\varepsilon$ cannot be bounded with respect to $t$.
\end{remarque}

Before we prove this result, we need a few technical lemmas:

\subsection{Technical preliminary results}

\begin{lemme}
  \label{lm:technique} Under the hypotheses of Proposition~\ref{pr:cv},
  for any $\varepsilon\in(0,1),$ we 
  have the following estimate:
  \begin{equation}
    \label{eq:trace}
    \int_{\Gamma_\varepsilon^0} u_\varepsilon^2 +\int_{\Gamma_\varepsilon^1} u_\varepsilon^2 + \varepsilon \int_{\Gamma_\varepsilon^\alpha}
    u_\varepsilon^2 \leq C\left( \int_{\Omega_\varepsilon} u_\varepsilon^2 +
      \int_{\Omega_\varepsilon} |\nabla u_\varepsilon|^2\right),
  \end{equation}
where $C>0$ does not depend on $\varepsilon.$
\end{lemme}
\noindent{\bf Proof:}
We use the same rescaling as in Section~\ref{sec:scale}, and define
\begin{equation}\label{eq:uv}
v(x,y) = u_\varepsilon\left(x,\varepsilon y \right).
\end{equation}
Then, the fact that
$u_\varepsilon\in H^1(\Omega_\varepsilon)$ implies that $v\in H^1(\Omega_1)$. On
this fixed domain, we can apply standard trace theorems
\cite{ding,gagliardo,lions-magenes}, which imply that there is a
constant $C>0$ depending on $\alpha$ only, such that
$$\int_{\partial \Omega_1} v^2 \leq C\left(\int_{\Omega_1} v^2+
  \int_{\Omega_1} |\nabla v|^2\right).$$
Inserting (\ref{eq:uv}) into this equation, we find that
\begin{multline*}
\int_{\Gamma^0_1} u_\varepsilon(x,\varepsilon y)^2 + \int_{\Gamma^1_1}
u_\varepsilon(x,\varepsilon y)^2 + \int_{\Gamma^\alpha_1} u_\varepsilon(x,\varepsilon y)^2 \leq
C \int_{\Omega_1} u_\varepsilon(x,\varepsilon y)^2\\
+ C\int_{\Omega_1}
\left(\frac{\partial u_\varepsilon}{\partial x}(x,\varepsilon  y) \right)^2 +
\varepsilon^2 \left(\frac{\partial u_\varepsilon}{\partial y}(x,\varepsilon  y)
\right)^2 
\end{multline*}
Hence, changing variables in these integrals, we have
$$\frac 1 \varepsilon \int_{\Gamma^0_\varepsilon} u_\varepsilon^2 + \frac 1
\varepsilon \int_{\Gamma^1_\varepsilon}
u_\varepsilon^2 + \int_{\Gamma^\alpha_\varepsilon} u_\varepsilon^2 \leq
C\int_{\Omega_\varepsilon} \frac 1 \varepsilon u_\varepsilon^2 + 
\frac 1 \varepsilon \left(\frac{\partial u_\varepsilon}{\partial x} \right)^2 +
\varepsilon \left(\frac{\partial u_\varepsilon}{\partial y}
\right)^2,$$
which proves the result.
\hfill$\Box$
\begin{remarque}
  In the above proof, we did not use the fact that $u_\varepsilon$
  satisfies (\ref{eq:simple}). Hence, the result of
  Lemma~\ref{lm:technique} is valid for any $u_\varepsilon \in
  H^1(\Omega_\varepsilon)$.
\end{remarque}
\begin{lemme}\label{lm:convergence}
  Under the hypotheses of Proposition~\ref{pr:cv}, there exists $u\in
  L^2([0,T]\times \Omega)$ such that $\nabla u \in L^2([0,T]\times
  \Omega^+)$ and $\nabla u \in L^2([0,T]\times
  \Omega^-)$, and the following convergences hold, up to extracting a
  subsequence:
  \begin{eqnarray}
    \label{eq:cv_u}
    u_\varepsilon  &\displaystyle \mathop{\longrightharpoonup}_{\varepsilon\to 0} & u
     \text{ in } L^2\left([0,T]\times \Omega\right), \\
    \label{eq:cv_grad_u}
    \nabla u_\varepsilon  &\displaystyle \mathop{\longrightharpoonup}_{\varepsilon\to
      0} & \nabla u  \text{ in } L^2\left([0,T]\times \Omega^+\right),
  \end{eqnarray}
\end{lemme}
\noindent{\bf Proof:} 
First note that the function $u_\varepsilon(t=0)$
 converges weakly to $u(t=0)$. Hence,
  \begin{equation}\label{eq:borne_0}
    \|u_\varepsilon(t=0)\|_{L^2(\Omega)} \leq C,
  \end{equation}
for some constant $C>0$ independent of $\varepsilon$. 
We consider (\ref{eq:simple}), multiply it by
$u_\varepsilon$, and integrate it over $\Omega_\varepsilon$ :
$$\frac12 \frac{d}{dt} \int_{\Omega_\varepsilon} u_\varepsilon^2(x,y,t)dxdy
-\int_{\Omega_{\varepsilon}} \Delta u_\varepsilon(x,y,t) u_\varepsilon(x,y,t)dxdy = 0.$$
Using an integration by parts and the boundary conditions in
(\ref{eq:simple}), we have
\begin{multline*}
-\int_{\Omega_\varepsilon} \Delta u_\varepsilon(x,y,t) u_\varepsilon(x,y,t)dxdy =
\int_{\Omega_\varepsilon} |\nabla u_\varepsilon(x,y,t)|^2dxdy -
\int_{\partial \Omega_\varepsilon} u_\varepsilon\partial_n
u_\varepsilon \\
= \int_{\Omega_\varepsilon} |\nabla u_\varepsilon(x,y,t)|^2dxdy -
\int_{\Gamma^1_\varepsilon} u_\varepsilon(x,y,t)dy
- \alpha\varepsilon \int_{\Gamma^\alpha_\varepsilon} u_\varepsilon(x,y,t)dx.
\end{multline*}
Thus, applying Cauchy-Schwarz inequality, 
\begin{multline*}
\frac12\frac{d}{dt} \left(\int_{\Omega_\varepsilon}
u_\varepsilon^2(x,y,t)dx\right)+\int_{\Omega_\varepsilon} |\nabla
u_\varepsilon(x,y,t)|^2dx  = \int_{\Gamma^1_\varepsilon} u_\varepsilon
+ \alpha\varepsilon \int_{\Gamma^\alpha_\varepsilon} u_\varepsilon \\
\leq \sqrt{\varepsilon(1-\alpha)}\left(\int_{\Gamma^1_\varepsilon} u_\varepsilon^2\right)^{1/2}
+ \alpha\varepsilon \left(\int_{\Gamma^\alpha_\varepsilon} u_\varepsilon^2\right)^{1/2}
\end{multline*}
We then apply Lemma~\ref{lm:technique}, finding
\begin{multline}
\label{eq:borne_u_0}
\frac12\frac{d}{dt} \left(\int_{\Omega_\varepsilon}
u_\varepsilon^2(x,y,t)dx\right)+\int_{\Omega_\varepsilon} |\nabla
u_\varepsilon(x,y,t)|^2dx 
\leq C\sqrt{\varepsilon(1-\alpha)}\left(\int_{\Omega_\varepsilon}
  |\nabla u_\varepsilon|^2\right)^{1/2} \\
+ C\alpha\sqrt\varepsilon \left(\int_{\Omega_\varepsilon} |\nabla u_\varepsilon|^2\right)^{1/2},
\end{multline}
for some constant $C$ depending only on $\alpha$. As a consequence,
there exists a constant $C$ (possibly different from the preceding one),
for which we have
$$\frac{d}{dt} \left(\int_{\Omega_\varepsilon}
u_\varepsilon^2(x,y,t)dx\right) \leq C\varepsilon.$$
Integrating this equation with respect to time, we thus have
$$\int_{\Omega_\varepsilon}
u_\varepsilon^2(x,y,t)dxdy \leq \int_{\Omega_\varepsilon}
u_\varepsilon^2(x,y,0)dxdy + C\varepsilon t.$$
Then, we split $\Omega$ into $1/\varepsilon$ domains of size
$\varepsilon$ in the direction $y$, and apply this inequality to each of
these domains. Since $u_\varepsilon = 0$ outside $\Omega_\varepsilon$,
this immediately implies
$$\int_{\Omega}
u_\varepsilon^2(x,y,t)dxdy \leq \int_{\Omega}
u_\varepsilon^2(x,y,0)dxdy + Ct \leq C(1+t),$$
for some constant $C>0$. The last inequality uses (\ref{eq:borne_0}). We
next integrate with respect to $t$, finding that the sequence
$u_\varepsilon$ is bounded independently of $\varepsilon$ in
$L^2(\Omega\times[0,T])$. Hence, up to extracting a subsequence, it
converges to some $u\in L^2(\Omega\times[0,T])$.

Next, going back to (\ref{eq:borne_u_0}), and integrating with respect
to time, we have
\begin{equation}\label{eq:borne_u_3}
\int_0^T \int_{\Omega_\varepsilon} |\nabla u_\varepsilon|^2(x,y,t)dxdydt \leq
C\varepsilon T + \frac12 \int_{\Omega_\varepsilon}u_\varepsilon^2(x,y,0)dxdy.
\end{equation}
This immediately implies that 
$$\int_0^T\int_{\Omega^+} |\nabla u_\varepsilon|^2(x,y,t)dxdydt \leq CT
+\frac12 \int_{\Omega} u_\varepsilon^2(x,y,0)dxdy.$$
Hence, $u_\varepsilon$ is bounded in $L^2([0,T],
H^1(\Omega^+))$. Extracting a subsequence if necessary, we thus have
(\ref{eq:cv_grad_u}).\hfill$\Box$

\begin{lemme}\label{lm:maj}
  Under the hypotheses of Proposition~\ref{pr:cv}, there exists a
  constant $C>0$ independent on $\varepsilon$ and $T$ such that:
  \begin{equation}
    \label{eq:borne_fond}
    \| u_\varepsilon\|_{L^2([0,T]\times \Gamma_\varepsilon^\alpha)} \leq C(T+1)
  \end{equation}
up to the extraction of a subsequence.
\end{lemme}
\noindent{\bf Proof:} We go back to (\ref{eq:borne_u_3}), which implies
that 
$$\int_0^T \int_{\Omega_\varepsilon} |\nabla u_\varepsilon|^2 \leq
C(T+1)\varepsilon,$$
where $C$ does not depend on $\varepsilon$ nor on $T$. Indeed,
$u_\varepsilon(t=0)$ satisfies (\ref{eq:borne_0}), and is $\varepsilon$-periodic with respect to $y$. Hence, 
$$\int_{\Omega_\varepsilon} u_\varepsilon^2 (x,y,0)dxdy \leq
C\varepsilon.$$
Here again, we use the scaling (\ref{eq:uv}), namely
$$v_\varepsilon(x,y,t) = u_\varepsilon(x,\varepsilon y,t),$$
and find that 
$$\int_0^T \int_{\Omega_1} \left(\partial_x v_\varepsilon\right)^2 +
\frac 1 {\varepsilon^2} \left(\partial_y v_\varepsilon\right)^2 \leq
C(1+T).$$
In particular, $v_\varepsilon$ is bounded in $L^2([0,T],
H^1(\Omega_1)).$ Using trace theorems \cite{ding,gagliardo,lions-magenes}, we infer that
$v_\varepsilon$ is bounded in $L^2([0,T], H^{1/2}(\Gamma_1^\alpha)),$
where
$$\Gamma_1^\alpha = \left\{ \left(x,\frac\alpha 2\right), \
  -1<x<0\right\}\bigcup \left\{ \left(x,-\frac\alpha 2\right), \
  -1<x<0\right\} .$$ 
In particular, $v_\varepsilon$ is bounded in $L^2([0,T]\times
\Gamma_1^\alpha)$. We finally point out that 
$$\|u_\varepsilon\|_{L^2([0,T]\times
\Gamma_\varepsilon^\alpha)} = \|v_\varepsilon\|_{L^2([0,T]\times
\Gamma_1^\alpha)},$$
which  completes the proof. \hfill$\Box$

\begin{lemme}
  \label{lm:cv_u_trace}
Under the hypotheses of Proposition~\ref{pr:cv}, we have the following
convergences, up to the extraction of a subsequence:
\begin{eqnarray}
    \label{eq:cv_u_bord_gauche}
    u_\varepsilon  &\displaystyle \mathop{\longrightharpoonup}_{\varepsilon\to
      0} & u  \text{ in } L^2\left([0,T]\times \Gamma^1 \right),
    \\
    \label{eq:cv_u_bord_droit}
    u_\varepsilon  &\displaystyle \mathop{\longrightarrow}_{\varepsilon\to
      0} & u  \text{ in } L^2\left([0,T]\times \Gamma^0 \right),
\end{eqnarray}
where $u$ is defined in Lemma~\ref{lm:convergence}.
Moreover, for any $\delta \in (0,1),$ we have the following
convergences, up to the extraction of a subsequence:
\begin{eqnarray}
    \label{eq:cv_u_bord_gauche_delta}
    u_\varepsilon  &\displaystyle \mathop{\longrightharpoonup}_{\varepsilon\to
      0} & u  \text{ in } L^2\left([0,T]\times \Omega\cap \{x=-\delta\} \right),
    \\
    \label{eq:cv_u_bord_droit_delta}
    u_\varepsilon  &\displaystyle \mathop{\longrightharpoonup}_{\varepsilon\to
      0} & u  \text{ in } L^2\left([0,T]\times \Omega\cap\{x=\delta\}
    \right),
  \end{eqnarray}
where $u$ is defined in Lemma~\ref{lm:convergence}. 
\end{lemme}
\noindent{\bf Proof:}
We already know that $u_\varepsilon$ converges weakly to $u$ in
$L^2([0,T],H^1(\Omega^+)).$ Using trace theorems, we infer that we have
weak convergence in $L^2([0,T], H^{1/2}(\Gamma^0)),$ hence strong
convergence in $L^2([0,T], L^2(\Gamma^0)) = L^2([0,T]\times \Gamma^0).$
This proves (\ref{eq:cv_u_bord_droit}).

Next, we prove (\ref{eq:cv_u_bord_gauche}). In view of (\ref{eq:trace}),
we already know that
$u_\varepsilon$ is bounded in $L^2([0,T]\times \Gamma^1)$, thus it
converges weakly, up to extracting a subsequence, to some limit. We are
now going to prove that this limit is $u$.
For this purpose, we use
here again the scaling (\ref{eq:uv}), namely
$$v_\varepsilon(x,y,t) = u_\varepsilon(x,\varepsilon y,t).$$
We use (\ref{eq:borne_u_3}), which implies that $v_\varepsilon$ is
bounded in $L^2([0,T], H^1(\Omega_1)).$ Hence, we have the following
convergence:
$$v_\varepsilon \mathop{\longrightharpoonup}_{\varepsilon \to 0} v \text
{ in } L^2([0,T],H^1(\Omega_1)),$$ for some $v\in
L^2([0,T],H^1(\Omega_1))$. Using the link between $u_\varepsilon$ and
$v_\varepsilon$, one easily proves using (\ref{eq:cv_u}), that 
$$u(x,y,t) = \int_{-\alpha/2}^{\alpha/2} v(x,z,t)dz.$$
Now, using trace theorems \cite{ding,gagliardo,lions-magenes}, we also have weak convergence of
$v_\varepsilon$ to $v$ in $L^2([0,T]\times H^{1/2}(\Gamma^1_1)),$ hence
in $L^2([0,T]\times\Gamma_1^1).$ Now, let $\varphi$ be a test function
in $C^\infty([0,T]\times \Gamma^1)$, and let us compute the integral of
$u_\varepsilon \varphi$ on $[0,T]\times\Gamma^1)$:
\begin{eqnarray*}
\int_0^T\int_{\Gamma^1} u_\varepsilon \varphi &=& \sum_{k\in \ZZ,
    \varepsilon|k|<1/2} \int_0^T \int_{\Gamma_\varepsilon^1 + \varepsilon k e_2}
  u_\varepsilon \varphi \\
&=& \sum_{k\in \ZZ,
    \varepsilon|k|<1/2} \int_0^T \int_{-\varepsilon\alpha/2}^{\varepsilon\alpha/2}
  u_\varepsilon(-1,y,t) \varphi(-1,y+k\varepsilon,t) dydt \\
&=&\sum_{k\in \ZZ,
    \varepsilon|k|<1/2} \int_0^T \int_{-\varepsilon\alpha/2}^{\varepsilon\alpha/2}
  v_\varepsilon\left(-1,\frac y\varepsilon,t\right)
  \varphi(-1,y+k\varepsilon,t) dydt \\
&=&\sum_{k\in \ZZ,
    \varepsilon|k|<1/2} \varepsilon \int_0^T \int_{-\alpha/2}^{\alpha/2}
  v_\varepsilon\left(-1,z,t\right) \varphi(-1,\varepsilon z+k\varepsilon,t) dzdt, 
\end{eqnarray*}
where we have used the fact that $u_\varepsilon$ is extended by $0$
outside $\Omega_\varepsilon$, the fact that $u_\varepsilon$ is
$\varepsilon$-periodic in $y$, and the link between $u_\varepsilon$ and
$v_\varepsilon$. Since $\varphi$ is smooth, one easily proves that 
$$\sum_{k\in \ZZ, \varepsilon|k|<1/2} \varepsilon \varphi(-1,\varepsilon
z + k\varepsilon,t) = \int_{-1/2}^{1/2} \varphi(-1,y,t)dy + O(\varepsilon),$$
where the remainder does not depend on $z$ nor on $v_\varepsilon$. Hence, since
$v_\varepsilon$ converges to $v$
in $L^2([0,T]\times \Gamma^1_1)$, we infer 
\begin{multline*}
\lim_{\varepsilon\to 0} \int_0^T\int_{\Gamma^1} u_\varepsilon\varphi =
\int_0^T \left(\int_{-\alpha/2}^{\alpha/2} v(-1,z,t)dz \int_{-1/2}^{1/2}
\varphi(-1,y,t) dy\right) dt  \\= \int_0^T \int_{-1/2}^{1/2}
u(-1,y,t)\varphi(-1,y,t) dydt.
\end{multline*}
We thus have proved (\ref{eq:cv_u_bord_gauche}). 

The convergence (\ref{eq:cv_u_bord_gauche_delta}) follows exactly the
same pattern. The proof of (\ref{eq:cv_u_bord_droit_delta}) is a direct
consequence of (\ref{eq:cv_grad_u}) and of trace theorems \cite{ding,gagliardo,lions-magenes}.\hfill$\Box$

\begin{lemme}
  \label{lm:cv-gamma-beta}
Under the hypotheses of Proposition~\ref{pr:cv}, we have the following
convergence, up to extraction of a subsequence:
\begin{equation}
  \label{eq:cv_u_bord_beta}
  \frac 1 \varepsilon\int_0^T \int_{\Gamma_\varepsilon^\beta} u_\varepsilon
  \mathop{\longrightarrow}_{\varepsilon\to 0} \alpha\int_0^T \int_{\Gamma^\beta} u.
\end{equation}
The function $u$ is defined in Lemma~\ref{lm:convergence}, and in the
right-hand side, $u$ is the trace on $\Gamma_\beta$ of $u_{|\Omega^+}$.
\end{lemme}
\noindent{\bf Proof:} We already know that $u_\varepsilon$ is bounded in $L^2([0,T],
H^1(\Omega^+))$. Hence, using trace theorems \cite{ding,gagliardo,lions-magenes},
$u_\varepsilon$ is bounded in $L^2([0,T],H^{1/2}(\{x=0\}\cap \Omega)) = L
^2([0,T], H^{1/2}(\Gamma^\beta)).$
Hence, up to extracting a subsequence, $u_\varepsilon$ converges in
$L^2([0,T]\times \Gamma^\beta).$ Since the trace operator is continuous,
its limit must be the trace of $u_{|\Omega^+}$. Now, since
$u_\varepsilon$ is $\varepsilon$-periodic in $y$, we have 
$$ \frac 1 \varepsilon\int_0^T\int_{\Gamma_\varepsilon^\beta} u_\varepsilon =
\sum_{k\in\ZZ, \varepsilon|k|<1/2} \int_0^T\int_{\Gamma_\varepsilon^\beta +
  k\varepsilon e_2} u_\varepsilon = \int_0^T\int_{\Gamma^\beta} u_\varepsilon
\sum_{k\in\ZZ, \varepsilon|k|<1/2} \11_{|y-k\varepsilon|<\alpha/2}.$$
In this integral, we have strong convergence on $u_\varepsilon$, whereas 
$$\sum_{k\in\ZZ, \varepsilon|k|<1/2} \11_{|y-k\varepsilon|<\alpha/2}
\mathop{\longrightharpoonup}_{\varepsilon\to 0} \alpha \text{ in }
L^2([0,T]\times \Gamma^\beta).$$
Hence, we may pass to the limit and obtain (\ref{eq:cv_u_bord_beta}).
\hfill$\Box$

\subsection{Proof of Proposition~\ref{pr:cv}}

We are now in position to give the

\noindent{\bf Proof of Proposition~\ref{pr:cv}:} We first apply
Lemmas~\ref{lm:convergence}, \ref{lm:maj} and \ref{lm:cv_u_trace},
getting the convergences
(\ref{eq:cv_u}), (\ref{eq:cv_grad_u}), (\ref{eq:borne_fond}) and (\ref{eq:cv_u_bord_gauche}).
We next prove that this limit $u$ is a solution to
(\ref{eq:hmodel.1})-(\ref{eq:hmodel.2})-(\ref{eq:htransmission}). For
this purpose, we define
\begin{equation}
  \label{eq:Omega_eps_-}
  \Omega_\varepsilon^- = \Omega_\varepsilon \cap \left\{ x<0\right\},
\end{equation}
and
\begin{equation}
  \label{eq:Omega_eps_+}
  \Omega_\varepsilon^+ = \Omega_\varepsilon \cap \left\{ x>0\right\}.
\end{equation}
We assume that $\varphi \in C^\infty(\Omega\times [0,T]),$ such that
$$\forall (x,y)\in \Omega, \ \varphi(x,y,T) = 0, \quad\text{and} \quad
\forall y\in \left(-\frac12,\frac12\right), \forall t\in [0,T], \
  \varphi(0,y,t) = \partial_x \varphi(0,y,t)= 0.$$ We multiply the equation
satisfied by $u_\varepsilon$ by $\varphi$, and integrate over
$\Omega_\varepsilon^+\times [0,T]$:
\begin{equation}
\label{eq:ipp}
\int_0^T \int_{\Omega_\varepsilon^+} \partial_t u_\varepsilon \varphi dxdydt-
\int_0^T\int_{\Omega_\varepsilon^+} \Delta u_\varepsilon\varphi dxdydt =0.
\end{equation}
In the first term, we integrate by parts with respect to $t$:
\begin{eqnarray}
\label{eq:ipp_1}
\int_0^T \int_{\Omega_\varepsilon^+} \partial_t u_\varepsilon \varphi
dxdydt &=& -\int_0^T \int_{\Omega_\varepsilon^+} u_\varepsilon  \partial_t\varphi
dxdydt \nonumber \\
&&- \int_{\Omega_\varepsilon^+} u_\varepsilon(x,y,0)
\varphi(x,y,0) dxdy 
\end{eqnarray}
Repeating this argument in each set $\Omega_\varepsilon^+ + k
\varepsilon e_2$, where $k\in \ZZ$ and $\varepsilon|k|<1/2$, and using
the fact that the union of these sets is $\Omega^+$ (if
$\varepsilon^{-1}$ is an integer, which we may assume), we have
\begin{eqnarray}
\label{eq:ipp_1.1}
\int_0^T \int_{\Omega^+} \partial_t u_\varepsilon \varphi
dxdydt &=& -\int_0^T \int_{\Omega^+} u_\varepsilon  \partial_t\varphi
dxdydt \nonumber \\
&&- \int_{\Omega^+} u_\varepsilon(x,y,0)
\varphi(x,y,0) dxdy
\end{eqnarray}
The weak convergence of $u_\varepsilon$ in $L^2(\Omega\times [0,T])$ and
of $u_\varepsilon(t=0)$ in $L^2(\Omega)$ allows to pass to the limit in
the right-hand side, finding
\begin{multline}
\label{eq:cv_u_dt_Omega_+}
  \lim_{\varepsilon\to 0} \sum_{k\in \ZZ, \varepsilon|k|<\frac12}
  \int_0^T \int_{\Omega_\varepsilon^+ + k\varepsilon e_2} \partial_t u_\varepsilon \varphi
dxdydt = \\-\int_0^T \int_{\Omega^+} u  \partial_t\varphi
dxdydt - \int_{\Omega^+} u(x,y,0)
\varphi(x,y,0) dxdy 
\end{multline}
The same argument allows to prove this convergence in $\Omega^-$ (recall
that we have extended $u_\varepsilon$ by $0$ outside $\Omega_\varepsilon$): 
\begin{equation}
\label{eq:cv_u_dt_Omega_-}
  \lim_{\varepsilon\to 0} \int_0^T \int_{\Omega_\varepsilon^-} \partial_t u_\varepsilon \varphi
dxdydt = -\int_0^T \int_{\Omega^-} u  \partial_t\varphi
dxdydt - \int_{\Omega^-} u(x,y,0)
\varphi(x,y,0) dxdy 
\end{equation}
We next deal with the second term in (\ref{eq:ipp}), in which we
integrate by parts with respect to $(x,y)$. Using the boundary
conditions we have on $u_\varepsilon$, we infer
\begin{multline}
 \label{eq:ipp_0.0}
  \int_0^T\int_{\Omega_\varepsilon^+} \Delta u_\varepsilon\varphi dxdydt =
  \int_0^T \int_{\partial \Omega_\varepsilon^+} \partial_n u_\varepsilon
  \varphi - \int_0^T \int_{\Omega_\varepsilon^+}  \nabla u_\varepsilon\cdot
  \nabla \varphi \\
=- \int_0^T \int_{\Omega_\varepsilon^+}  \nabla u_\varepsilon\cdot
  \nabla \varphi.
\end{multline}
Using (\ref{eq:cv_grad_u}), we thus have
$$\lim_{\varepsilon\to 0} \int_0^T\int_{\Omega_\varepsilon^+} \Delta
u_\varepsilon\varphi dxdydt = - \int_0^T \int_{\Omega^+}  \nabla
u_\varepsilon\cdot \nabla \varphi.$$
This and (\ref{eq:cv_u_dt_Omega_+}) implies 
\begin{equation}
  \label{eq:pb_Omega_+}
   -\int_0^T \int_{\Omega^+} u  \partial_t\varphi
dxdydt - \int_{\Omega^+} u(x,y,0)
\varphi(x,y,0) dxdy= - \int_0^T \int_{\Omega^+}  \nabla
u_\cdot \nabla \varphi,
\end{equation}
which is the weak formulation of (\ref{eq:hmodel.2}).

\medskip

We are now going to apply the same strategy on the set
$\Omega_\varepsilon^-$, but the situation is more delicate
here. Integrating over $\Omega_\varepsilon^-$ instead of
$\Omega_\varepsilon^+$, we already know that (\ref{eq:cv_u_dt_Omega_-})
holds. Moreover, we have
\begin{multline}
 \label{eq:ipp_0}
  \int_0^T\int_{\Omega_\varepsilon^-} \Delta u_\varepsilon\varphi dxdydt =
  \int_0^T \int_{\partial \Omega_\varepsilon^-} \partial_n u_\varepsilon
  \varphi - \int_0^T \int_{\partial \Omega_\varepsilon^-}  u_\varepsilon
  \partial_n \varphi + \int_0^T \int_{\Omega_\varepsilon^-} u_\varepsilon
  \Delta \varphi \\
= \int_0^T \int_{\Gamma_\varepsilon^1} \varphi + \int_0^T
(\alpha-\beta)\frac{\varepsilon}2 \int_{\Gamma_\varepsilon^\alpha} \varphi -
\int_0^T \int_{\partial \Omega_\varepsilon^-} u_\varepsilon\partial_n\varphi + \int_0^T
\int_{\Omega_\varepsilon^-} u_\varepsilon \Delta \varphi.
\end{multline}
In order to recover an integral over the domain $\Omega^-$, we repeat this
operation for the domain $\Omega_\varepsilon^- + \varepsilon k,$ for $k\in
\ZZ$ and $\varepsilon |k| \leq 1/2.$ We thus have (\ref{eq:ipp_0}) for
this domain. Summing these equalities with respect to $k$, and setting
\begin{equation}
\label{eq:Omega_tilde}
\tilde\Omega_\varepsilon^- = \bigcup_{k\in \ZZ, \ \varepsilon |k| < 1/2}
\left(\Omega_\varepsilon^- + \varepsilon ke_2\right),
\end{equation}
we find
\begin{multline}
 \label{eq:ipp_3}
  \int_0^T\int_{\tilde \Omega_\varepsilon^-} \Delta u_\varepsilon\varphi dxdydt 
= \sum_{k\in \ZZ, \ \varepsilon |k| < 1/2} \int_0^T
\int_{\Gamma_\varepsilon^1+\varepsilon k e_2} \varphi 
\\+ \sum_{k\in \ZZ, \ \varepsilon |k| < 1/2}\int_0^T
(\alpha-\beta)\frac{\varepsilon}2 \int_{\Gamma_\varepsilon^\alpha+\varepsilon ke_2}
\varphi -
\sum_{k\in \ZZ, \ \varepsilon |k| < 1/2}\int_0^T \int_{\partial
  (\Omega_\varepsilon^-+\varepsilon ke_2)} u_\varepsilon\partial_n\varphi \\+ \int_0^T
\int_{\Omega^-} u_\varepsilon \Delta \varphi.
\end{multline}
We deal with each term of the right-hand side of (\ref{eq:ipp_3}) separately: for the first
term, we have
\begin{equation}
  \label{eq:ipp_4}
 \sum_{k\in \ZZ, \ \varepsilon |k| < 1/2} \int_0^T
\int_{\Gamma_\varepsilon^1+\varepsilon ke_2} \varphi
\mathop{\longrightarrow}_{\varepsilon \to 0} (1-\alpha)\int_0^T \int_{\Gamma^1}\varphi.
\end{equation}
Next, the second term of the right-hand side of (\ref{eq:ipp_3}) is interpreted as a Riemann
sum with respect to $k$ for the map $y\mapsto \int_{\Gamma^1_\varepsilon+ye_2}\varphi$. Since $\varphi$ is smooth, we
thus have
\begin{equation}
  \label{eq:ipp_5}
\sum_{k\in \ZZ, \ \varepsilon |k| < 1/2}\int_0^T
(\alpha-\beta)\frac{\varepsilon}2 \int_{\Gamma_\varepsilon^\alpha+\varepsilon ke_2}
\varphi \mathop{\longrightarrow}_{\varepsilon\to 0}
(\alpha-\beta)\int_0^T\int_{\Omega^-} \varphi.  
\end{equation}
Finally, we deal with the third term of the right-hand side of
(\ref{eq:ipp_3}). We first note that, due to (\ref{eq:Omega_tilde}),
this term is equal to 
\begin{multline}
\label{eq:ipp_6}
\sum_{k\in \ZZ, \ \varepsilon |k| < 1/2}\int_0^T \int_{\partial
  (\Omega^-_\varepsilon+\varepsilon ke_2)} u_\varepsilon\partial_n\varphi =
\int_0^T \int_{\partial \tilde\Omega^-_\varepsilon} u_\varepsilon \partial_n\varphi
\\
= -\sum_{\varepsilon |k| <1/2} \int_0^T
\int_{\Gamma_\varepsilon^1+\varepsilon k e_2} u_\varepsilon\partial_x \varphi +
\sum_{\varepsilon |k|<1/2} \int_0^T \int_{\Gamma_\varepsilon^\alpha +
  \varepsilon k e_2} u_\varepsilon \partial_n \varphi 
\end{multline}
Here again, each term of the right-hand side of (\ref{eq:ipp_6}) is
dealt will separately. The first term reads
$$-\sum_{\varepsilon |k| <1/2} \int_0^T
\int_{\Gamma_\varepsilon^1+\varepsilon k e_2} u_\varepsilon\partial_x\varphi = -\int_0^T
\int_{\Gamma^1}
u_\varepsilon\partial_x\varphi.$$
Hence, applying (\ref{eq:cv_u_bord_gauche}), we find
\begin{equation}
  \label{eq:ipp_7}
 - \sum_{\varepsilon |k| <1/2} \int_0^T
\int_{\Gamma_\varepsilon^1+\varepsilon k e_2} u_\varepsilon\partial_x\varphi
\mathop{\longrightarrow}_{\varepsilon\to 0} - \int_{\Gamma^1}
u\partial_x \varphi.
\end{equation}
Turning to the second term of the right-hand side of (\ref{eq:ipp_6}),
we note that the value of $u_\varepsilon$ on this boundary does not
depend on $k$ since $u_\varepsilon$ is periodic of period $\varepsilon$
with respect to $y$. Hence, we write
\begin{multline*}
\sum_{\varepsilon |k|<1/2} \int_0^T \int_{\Gamma_\varepsilon^\alpha +
  \varepsilon k e_2} u_\varepsilon \partial_n\varphi = \\ \sum_{\varepsilon |k|<1/2}
\int_0^T \int_{-1}^0 u_\varepsilon\left(x,-\frac{\alpha}2\varepsilon,t\right)
\partial_y\varphi\left(x,\varepsilon k -\frac{\alpha}2 \varepsilon,t\right)dxdt \\
- \sum_{\varepsilon |k|<1/2}
\int_0^T \int_{-1}^0 u_\varepsilon\left(x,\frac{\alpha}2\varepsilon,t\right)
\partial_y\varphi\left(x,\varepsilon k +\frac{\alpha}2 \varepsilon,t\right)dxdt.
\end{multline*}
We use a Taylor expansion of $\varphi$ with respect to the $y$ variable,
around the point $y = \varepsilon k$:
$$\partial_y\varphi\left(x,\varepsilon k \pm\frac{\alpha}2 \varepsilon,t\right) =
\partial_y\varphi\left(x,\varepsilon k \varepsilon,t\right) \pm
\varepsilon\frac{\alpha}2 \partial_y^2\varphi(x,\varepsilon k, t) +
O\left(\varepsilon^2\right).
$$
The remainder $O(\varepsilon^2)$ depends only on $\varphi$, and may be
chosen bounded independently of $x$ and $t$. Hence, we have
\begin{multline*}
\sum_{\varepsilon |k|<1/2} \int_0^T \int_{\Gamma_\varepsilon^\alpha +
  \varepsilon k e_2} u_\varepsilon \partial_y \varphi = O(\varepsilon) \\
 +\sum_{\varepsilon |k|<1/2}
\int_0^T \int_{-1}^0
\left(u_\varepsilon\left(x,-\frac{\alpha}2\varepsilon,t\right) -
  u_\varepsilon\left(x,\frac{\alpha}2\varepsilon,t\right) \right)
\partial_y \varphi\left(x,\varepsilon k,t\right)dxdt \\
- \varepsilon\frac{\alpha}2 \sum_{\varepsilon |k|<1/2}
\int_0^T \int_{-1}^0
\left(u_\varepsilon\left(x,-\frac{\alpha}2\varepsilon,t\right) +
  u_\varepsilon\left(x,\frac{\alpha}2\varepsilon,t\right) \right)
\partial_y^2\varphi\left(x,\varepsilon k,t\right)dxdt.
\end{multline*}
The last term is a Riemann sum with respect to $k$. The corresponding
integral is 
$$\int_{-1/2}^{1/2}\partial_y^2\varphi(x,y,t)dy = 0.$$ 
Hence, using Lemma~\ref{lm:maj}, one easily proves that this
last term converges to $0$. For the first term, we use the same
strategy, but we need to prove a little more: we write the sum over $k$
as follows:
\begin{multline*}
\sum_{\varepsilon |k|<1/2}\partial_y\varphi(x,\varepsilon k,t) =
\sum_{\varepsilon |k|<1/2} \left(\partial_y\varphi(x,\varepsilon k,t)
  -\frac 1 \varepsilon \int_{\varepsilon k - \varepsilon/2}^{\varepsilon k +
    \varepsilon/2}\partial_y \varphi(x,y,t)dy\right) \\
= \sum_{\varepsilon |k|<1/2}  \int_{\varepsilon k - \varepsilon/2}^{\varepsilon k +
    \varepsilon/2}\frac{\partial_y\varphi(x,\varepsilon k,t)
  -\partial_y \varphi(x,y,t)}\varepsilon dy \\
= \sum_{\varepsilon |k|<1/2}  \int_{\varepsilon k - \varepsilon/2}^{\varepsilon k +
    \varepsilon/2} (y-\varepsilon k) \partial_y^2 \varphi(x,\varepsilon
  k, t)dy \\
+ \varepsilon\sum_{\varepsilon |k|<1/2}  \int_{\varepsilon k - \varepsilon/2}^{\varepsilon k +
    \varepsilon/2} \frac12(y-\varepsilon k)^2 \partial_y^3 \varphi(x,\varepsilon
  k, t)dy + O(\varepsilon),
\end{multline*}
where the remainder $O(\varepsilon)$ is uniform with respect to $x$ and
$t$. The first term vanishes because the integrand is even with respect
to the variable $y-\varepsilon k$. The second term is easily shown to be
of order $\varepsilon$ by computing the integrals explicitly. Finally,
we thus have
\begin{equation}
  \label{eq:ipp_9}
  \sum_{\varepsilon |k|<1/2} \int_0^T \int_{\Gamma_\varepsilon^\alpha +
  \varepsilon k e_2} u_\varepsilon \partial_n\varphi
\mathop{\longrightarrow}_{\varepsilon \to 0} 0.
\end{equation}
Inserting (\ref{eq:ipp_7}) and (\ref{eq:ipp_9}) into (\ref{eq:ipp_6}),
we infer
\begin{equation}
  \label{eq:ipp_10}
 \sum_{k\in \ZZ, \ \varepsilon |k| < 1/2}\int_0^T \int_{\partial
  (\Omega_\varepsilon+\varepsilon ke_2)} u_\varepsilon\partial_n\varphi
\mathop{\longrightarrow}_{\varepsilon\to 0} -\int_{\Gamma^1} u\partial_x \varphi.
\end{equation}
Then, we collect (\ref{eq:ipp_4}), (\ref{eq:ipp_5}), (\ref{eq:ipp_10}),
and insert them into (\ref{eq:ipp_3}). Hence,
$$  \int_0^T\int_{\tilde \Omega_\varepsilon^-} \Delta
u_\varepsilon\varphi dxdydt \mathop{\longrightarrow}_{\varepsilon\to 0}
(1-\alpha) \int_0^T \int_{\Gamma^1}\varphi + (\alpha-\beta)\int_0^T
\int_{\Omega^-} \varphi -\int_{\Gamma^1} u\partial_x\varphi - \int_0^T
\int_{\Omega^-} u\Delta \varphi.$$
Finally, we use this convergence and (\ref{eq:cv_u_dt_Omega_-}), and
insert it into the equation 
$$\int_0^T \int_{\Omega_\varepsilon^-} \partial_t u_\varepsilon \varphi dxdydt-
\int_0^T\int_{\Omega_\varepsilon^-} \Delta u_\varepsilon\varphi dxdydt
=0,$$
finding
\begin{multline}
  \label{eq:pb_Omega_-}
   -\int_0^T \int_{\Omega^-} u  \partial_t\varphi
dxdydt - \int_{\Omega^-} u(x,y,0)
\varphi(x,y,0) dxdydt= \\
(1-\alpha) \int_0^T \int_{\Gamma^1}\varphi + (\alpha-\beta)\int_0^T
\int_{\Omega^-} \varphi + \int_0^T
\int_{\Omega^-} \nabla u\cdot\nabla \varphi,
\end{multline}
which is a weak formulation of (\ref{eq:hmodel.1}).

\medskip

To end the proof, we need to show that the transmission conditions
(\ref{eq:htransmission}) hold. For this purpose, we first point out
that, since $u_\varepsilon \longrightharpoonup u$ in $L^2([0,T]\times
\Omega)$ and since $u_\varepsilon$ is $\varepsilon$-periodic in $y$,
$u$ must be independent of $y$. Next, we define a test function
$\varphi$ which depends only on $t$, and has 
compact support in $(0,T)$. We multiply
the first line of (\ref{eq:simple}) by $\varphi$ and integrate over
$[0,T] \times \tilde \Omega_{\varepsilon,\delta},$ where
$$\tilde\Omega_{\varepsilon,\delta} = \tilde\Omega_\varepsilon \cap
\left\{ |x|< \delta \right\}.$$
(Recall that $\tilde \Omega_\varepsilon$ is defined by
(\ref{eq:Omega_tilde}).) Integrating by parts, we have
\begin{eqnarray*}
  0 &=& \int_0^T \int_{\tilde \Omega_{\varepsilon,\delta}}
  \left(\partial_t u_\varepsilon - \Delta u_\varepsilon\right)
  \varphi(t) dxdydt \\ 
&=& -\int_0^T \int_{\tilde \Omega_{\varepsilon,\delta}}
  u_\varepsilon \partial_t \varphi - \int_0^T \int_{\partial \left( \tilde
    \Omega_{\varepsilon,\delta}\right)} \partial_n
    u_\varepsilon \varphi 
  \end{eqnarray*}
Computing the boundary term, we get
\begin{eqnarray}\label{eq:ipp_11.1}
0&=& -\int_0^T \int_{\tilde \Omega_{\varepsilon,\delta}}
  u_\varepsilon \partial_t \varphi 
  -\int_0^T \int_{\Omega\cap \{x = \delta\}} \partial_x
  u_\varepsilon\varphi + \int_0^T \int_{\tilde \Omega_\varepsilon\cap
    \{x = -\delta\}}\partial_x u_\varepsilon \varphi \nonumber \\
&& - (\alpha-\beta)\frac\varepsilon 2\sum_{k\in \ZZ, \varepsilon|k|< 1/2}
  \int_0^T \int_{-\delta}^0 2\varphi(t)dxdt
  \nonumber \\
&&- \frac\beta\alpha\sum_{k\in \ZZ, \varepsilon|k|< 1/2}
  \int_0^T\int_{\Gamma_\varepsilon^\beta + k\varepsilon e_2}\varphi.
\end{eqnarray}
It is easy to pass to the limit in each of the above terms, except for
the second and third one. We thus deal with them separately: we multiply
the equation by $\varphi$ and
integrate on $\tilde\Omega_\varepsilon\cap \{x<-\delta\}$
instead of $\tilde\Omega_{\varepsilon,\delta}$. Integrating by parts, we
find
\begin{eqnarray*}
0 &=& -\int_0^T\int_{\tilde\Omega_\varepsilon\cap\{x<-\delta\}}
u_\varepsilon \partial_t\varphi
-\int_0^T\int_{\tilde\Omega_\varepsilon\cap\{x=-\delta\}} \partial_x
u_\varepsilon\varphi \\
&& -(\alpha-\beta) \frac\varepsilon 2 \sum_{k\in \ZZ, \varepsilon|k|< 1/2}\int_0^T \int_{-1}^{-\delta} 2\varphi(t)dxdt.
\end{eqnarray*}
Passing to the limit in the first and third term, we find that 
$$\lim_{\varepsilon\to 0} \int_0^T\int_{\tilde\Omega_\varepsilon\cap\{x=-\delta\}} \partial_x
u_\varepsilon\varphi = -\int_0^T\int_{\tilde\Omega_\varepsilon\cap\{x<-\delta\}}
u \partial_t\varphi + (\alpha-\beta) \int_0^T \int_{\Omega\cap
  \{x<-\delta\}} \varphi.$$
Integrating by parts in the right-hand side and using the equation
satisfied by $u$, namely (\ref{eq:hmodel.1}), we infer
\begin{equation}
  \label{eq:ipp_11}
  \lim_{\varepsilon\to 0} \int_0^T\int_{\tilde\Omega_\varepsilon\cap\{x=-\delta\}} \partial_x
u_\varepsilon\varphi = \int_0^T
\int_{\Omega\cap\{x=-\delta\}} \partial_x u \varphi.
\end{equation}
Likewise, integrating over $\tilde\Omega_\varepsilon \cap\{x>\delta\}$,
we have
\begin{equation}
  \label{eq:ipp_12}
  \lim_{\varepsilon\to 0} \int_0^T\int_{\tilde\Omega_\varepsilon\cap\{x=\delta\}} \partial_x
u_\varepsilon\varphi = \int_0^T
\int_{\Omega\cap\{x=\delta\}} \partial_x u \varphi.
\end{equation}
We insert (\ref{eq:ipp_11}) and (\ref{eq:ipp_12}) into
(\ref{eq:ipp_11.1}), and interpreting the last term of
(\ref{eq:ipp_11.1}) as Riemann sum, we get
\begin{multline*}
0 = -\int_0^T\int_{\Omega\cap \{|x|<\delta\}} u\partial_t\varphi
   - \int_0^T \int_{\Omega\cap
  \{x=\delta\}} \partial_x u \varphi + \int_0^T
\int_{\Omega\cap\{x=-\delta\}} \partial_x u \varphi \\- (\alpha-\beta)\delta
\int_0^T \varphi(t)dt - \beta\int_0^T
\int_{\Gamma^\beta} \varphi.$$
\end{multline*}
Letting $\delta\to 0$, and using the fact that $u$ does not depend on
$y$, we thus find
\begin{equation}
  \label{eq:transmission_du}
\partial_x u^- = \partial_x u^+ + \beta \quad \text{if } x=0.  
\end{equation}
Then, we repeat the same argument with $x\varphi(t)$ instead of
$\varphi(t)$, and we find
$$0 = -\int_0^T \int_{\Omega\cap \{|x|<\delta\}} xu_\varepsilon \partial_t
  \varphi + \int_0^T \int_{\partial\left(\tilde
      \Omega_{\varepsilon,\delta}\right)}
  u_\varepsilon\partial_n\left(x\varphi\right) - x\varphi\partial_n
  u_\varepsilon .
$$
Since $x\varphi$ does not depend on $y$, the first boundary term
contains only terms on the vertical boundaries $\{x=\delta\}$,
$\tilde\Omega_\varepsilon\cap\{x=-\delta\},$ and
$\Gamma_\varepsilon^\beta + k\varepsilon e_2$, for $k\in \ZZ,$
$\varepsilon|k| < 1/2.$ We thus have
\begin{eqnarray*}
  0 &=& -\int_0^T \int_{\Omega\cap \{|x|<\delta\}} xu_\varepsilon \partial_t
  \varphi \\ 
&&+ \int_0^T \int_{\Omega\cap\{x=\delta\}} u_\varepsilon\varphi -
   \int_0^T\int_{\tilde \Omega_{\varepsilon}\cap\{x=-\delta\}}
  u_\varepsilon\varphi - \sum_{k\in \ZZ, \varepsilon|k|<1/2}
   \int_0^T\int_{\Gamma_\varepsilon^\beta + \varepsilon k e_2}
  u_\varepsilon\varphi \\
&&  -\int_0^T \int_{\Omega\cap \{x = \delta\}} \partial_x
  u_\varepsilon x\varphi + \int_0^T \int_{\tilde \Omega_\varepsilon\cap
    \{x = -\delta\}}\partial_x u_\varepsilon x\varphi \\
&& - (\alpha-\beta)\frac\varepsilon 2\sum_{k\in \ZZ, \varepsilon|k|< 1/2}
  \int_0^T \int_{-\delta}^0 2x\varphi(t)dxdt
\end{eqnarray*}
Here again, we use (\ref{eq:cv_u_bord_gauche_delta}),
(\ref{eq:cv_u_bord_droit_delta}), (\ref{eq:ipp_11}) and
(\ref{eq:ipp_12}), to pass to the limit in the above equation,
finding
\begin{eqnarray*}
0 &=& -\int_0^T \int_{\Omega\cap\{|x|< \delta\}} xu\partial_t\varphi +  \int_0^T\int_{\Omega\cap\{x=\delta\}} u\varphi -
   \int_0^T\int_{\Omega\cap\{x=-\delta\}}
  u\varphi - \alpha \int_0^T\int_{\Gamma^\beta}
  u^+\varphi \\
&&- \int_0^T \int_{\Omega\cap \{x = \delta\}} \partial_x
  u_\varepsilon x\varphi + \int_0^T \int_{\tilde \Omega_\varepsilon\cap
    \{x = -\delta\}}\partial_x u_\varepsilon x\varphi -
  (\alpha-\beta)\frac{\delta^2}2 \int_0^T \varphi  
\end{eqnarray*}
Hence, letting $\delta \to 0$, we have
\begin{equation}
  \label{eq:transmission_u}
  u^+ - u^- - \alpha u^+ = 0 \quad \text{if } x=0.
\end{equation}
Collecting (\ref{eq:transmission_du}) and (\ref{eq:transmission_u}), we
find (\ref{eq:htransmission}).

\medskip

Thus, we have proved that, up to extracting a subsequence, the sequence
$u_\varepsilon$ converges to $u$, weakly in $L^2$, and that $u\in H^1(\Omega)$ is
solution to
(\ref{eq:hmodel.1})-(\ref{eq:hmodel.2})-(\ref{eq:htransmission}). Next,
we point out that, according to Lemma~\ref{lm:existence_unicite_homog}, such a solution is
unique. Hence, the whole sequence converges, without any need to extract
a subsequence.\hfill$\Box$

\newpage

\section{A numerical illustration}

In this section, we illustrate the model MOSAIC on the crack geometry
given by figure~\ref{fig:dessin}, which is simple enough to be calculated directly. We can thus compare the solution of the diffusion problem (\ref{eq:simple}) that we will compute for different periods $\varepsilon\to 0$ with the homogenized approach, corresponding to (\ref{eq:hmodel.1})-(\ref{eq:hmodel.2}). To calculate the solution of the homogenized problem, we use two independent methods that we describe in details:
\begin{itemize}
\item{a weak formulation in the whole homogenized model $\Omega$.}
\item{a fixed-point approach connecting the 2 sub-domains $\Omega^+$ and $\Omega^-$.}
\end{itemize}
Computations presented here are performed using a $P^1$ finite element approximation on triangular meshes. It  has been implemented using the software FreeFem++ \cite{free07}.

\subsection{Direct simulation of the cracked domain}

\begin{figure}

        \centering
        \begin{subfigure}[b]{0.65\textwidth}
                \centering
                \includegraphics[width=\textwidth]{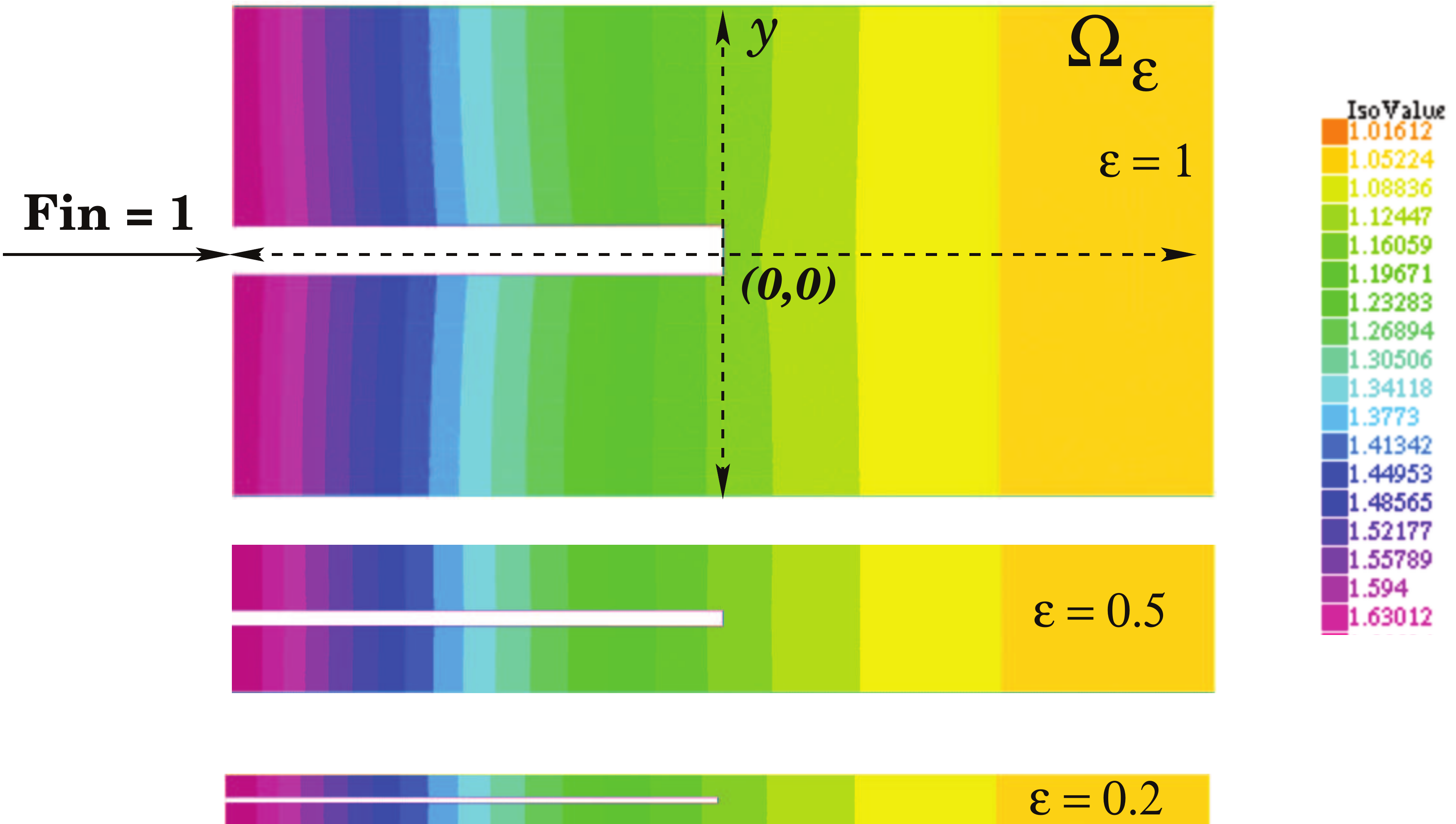}

                \caption{\textbf{$-\Delta u + \partial_t u =  0$ in $\Omega_\varepsilon$}.}
                \label{fig:direct}
        \end{subfigure}%
        ~ 
        \begin{subfigure}[b]{0.35\textwidth}
                \centering
                \includegraphics[width=\textwidth]{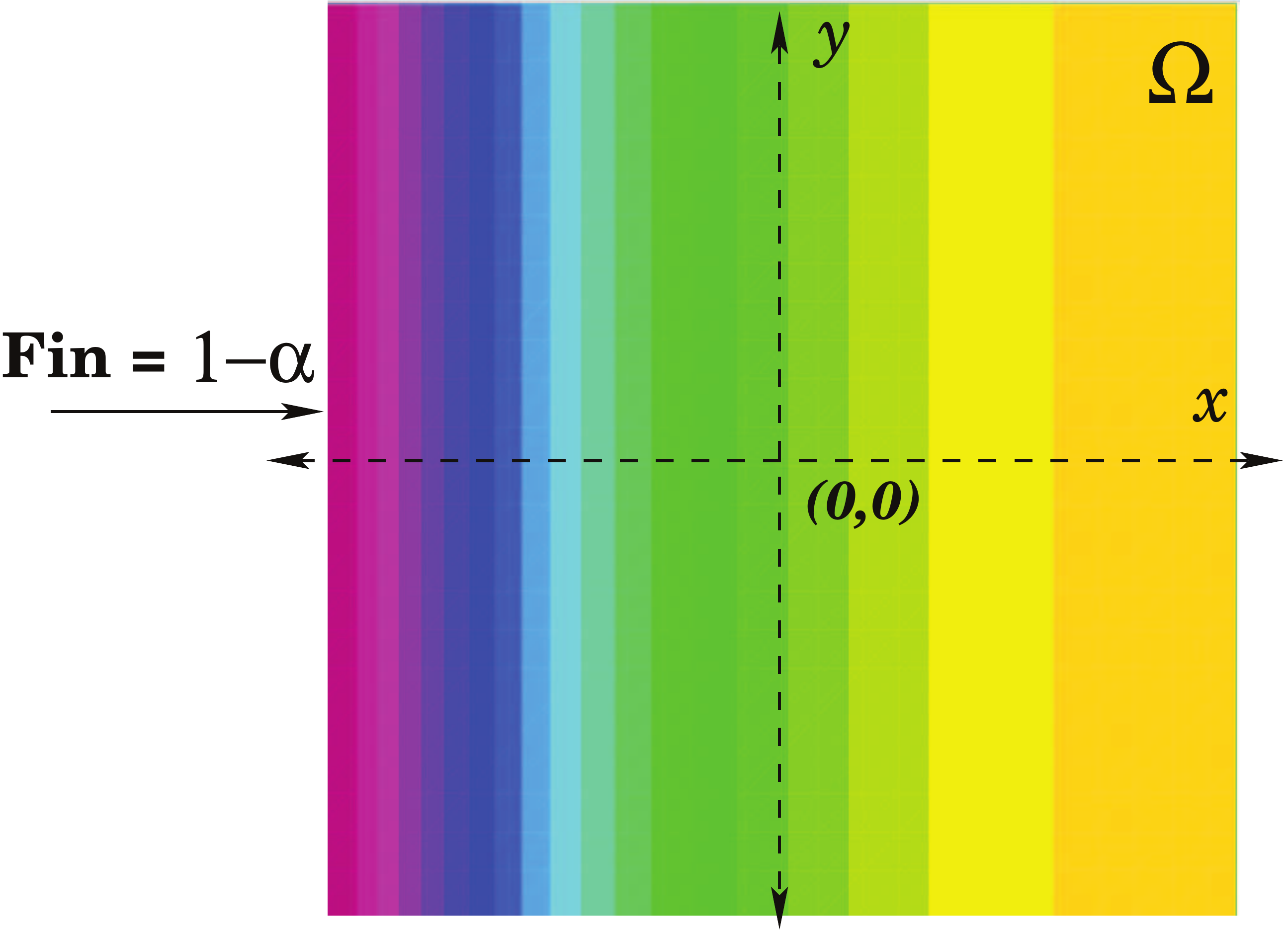}
                \caption{\textbf{$-\Delta u + \partial_t u =  \alpha\11_{\{x<0\}}$ in $\Omega$}.}
                \label{fig:approxh}

        \end{subfigure}
         \caption{\textbf{Direct calculation on the periodic cracked domain $\Omega_\varepsilon$ for different  periods $\varepsilon\to 0$}. \textit{The crack occupies the fraction $\alpha = 0.1$. The field solution of (\ref{eq:simple}) is plotted at time $t_1=0.5$. We compare the result with that of the homogenized approach approximated by equation (\ref{eq:hmodel}) where the singularity at $x=0$ is disregarded. This approximation makes sense here because $\alpha$ is rather small.}}
\label{fig:conv2}
\end{figure}

We start by solving the problem corresponding to the diffusion model
(\ref{eq:simple}) on a cracked domain $\Omega_\varepsilon$  whose shape
reproduces figure~\ref{fig:dessin}. We recall that the domain
$\Omega_\varepsilon$ is periodic of period $\varepsilon$ with respect to
the $y$ direction. The  width of the crack $\alpha$ is fixed: for
example, we take $\alpha = 0.1$ for the simulations illustrated on
figure ~\ref{fig:conv2}. In what follows, we always take $\beta=0$.
We represent the field $u(x,y,t_1)$ solution of
(\ref{eq:simple}) at a given time $t_1=0.5$ for different periods
$\varepsilon$.  The time $t_1$ is such that the normalized spatial
profile of the field has reached a stationary state. We carry out direct
calculations for the periods $\varepsilon = 1$, $\varepsilon = 0.5$,
$\varepsilon = 0.2$, $\varepsilon = 0.02$. We make sure the simulations
are converged with respect to the mesh size as well as the time
step. The fields calculated for each period $\varepsilon$ are shown on
figure~\ref{fig:direct}. We note that the convergence in $\varepsilon$ is
quite fast: for $\varepsilon\leq 0.2$, the field does not depend on the variable $y$ anymore and the limit in $\varepsilon$ appears to be reached. This is not true for $\varepsilon = 1$ where the field still depends on $y$, especially in the vicinity of the crack.
Besides, since $\alpha$ is rather small here, we make use of remark
\ref{rk:hmodel} and compare the field solution of the direct calculation
with that of the approximated homogenized model
(\ref{eq:hmodel}). Equation (\ref{eq:hmodel}) is solved in the domain
$\Omega$ where the crack is no more described in the geometry but its
effect on the diffusion process is modeled by a source term of the form
$\alpha\11_{\{x<0\}}$ in the right hand side of equation
(\ref{eq:hmodel}). The incoming flux applied on the left boundary
located at $x=-1$ is now $1-\alpha$ (see figure~\ref{fig:approxh}). It
is consistent with the fact that the indentation of width $\alpha$ is no
more described in the homogenized domain $\Omega$. A fraction of the
incoming flux is thus "converted" into a source term in the homogenized equation. This accounts for the name given to our approach: MOSAIC, as Model Of Sinks Averaging Inhomogeneous behavior of Cracked media.
Note that the homogenized model (\ref{eq:hmodel}) does not take the singularity at $x=0$ into account, which should play a significant role as $\alpha$ increases. We are now going to devise two approaches to treat this singularity in the homogenized model for any given $\alpha$.




\subsection{Weak formulation of the homogenized model}
\label{formulation-faible}
We compare the direct approach, considered as a reference point,
with the homogenized model established in
(\ref{eq:hmodel.1})-(\ref{eq:hmodel.2}) on the average domain $\Omega$
(i.e without any crack). This domain, which represents the limit of $\Omega_\varepsilon$ as $\varepsilon \to 0$, is invariant in the $y$ direction and its length is $2$ in the $x$ direction (same length as the cracked domain $\Omega_\varepsilon$).
With respect to the model (\ref{eq:hmodel.1})-(\ref{eq:hmodel.2}), the incoming flux applied on the left boundary located at $x=-1$ is now $1-\alpha$.


The first strategy to compute the solution of (\ref{eq:hmodel.1})-(\ref{eq:hmodel.2}) is to find an equivalent weak formulation of the problem and then solve it numerically using standard finite element methods. 

We recall the equations that have to be solved:
\begin{equation}
\label{eq:hmodel.1b}
\begin{cases}
  -\Delta u + \partial_t u =  \alpha-\beta &
  \text{in } \Omega^- = \{ -1<x<0\}, \\
  \partial_n u = 1- \alpha & \text{on } \{x=-1\}, \\
  u \text{ is }1-\text{periodic in } y.
\end{cases}
\end{equation}
\begin{equation}
\label{eq:hmodel.2b}
\begin{cases}
  -\Delta u + \partial_t u =  0 &
  \text{in } \Omega^+ = \{ 0<x<1\}, \\
  \partial_n u = 0 & \text{on } \{x=1\}, \\
  u \text{ is }1-\text{periodic in } y.
\end{cases}
\end{equation}
with the so-called transmission conditions at the bottom of the crack located at $x=0$:
\begin{equation}
  \label{eq:htransmissionb}
  u(x=0^-) = (1-\alpha)u(x=0^+), \quad \partial_x u(x=0^-)
  = \partial_x u(x=0^+) + \beta.
\end{equation}

Multiplying (\ref{eq:hmodel.1b}) by a test function $v$ and integrating by parts on the domain $\Omega^-$, we get:

$$\int_{\Omega^-} \partial_t uv -\Delta v u + \int_{\partial \Omega^-}
u\partial_n v - v\partial_n u = (\alpha - \beta) \int_{\Omega^-} v.$$

Using the boundary conditions on $\Omega^-$, this becomes:

\begin{multline}
  \label{eq:1}
  \int_{\Omega^-} \partial_t uv -\Delta v u + \int_{x=-1}
u\partial_n v - v(1-\alpha) + \int_{x=0} u(x=0^-)\partial_x v
-v \partial_x u(0^-)\\= (\alpha - \beta) \int_{\Omega^-} v.
\end{multline}

Then, applying the same processes (\ref{eq:hmodel.2b}) and integrating on $\Omega^+$, we obtain :
$$\int_{\Omega^+} \partial_t uv -\Delta v u + \int_{\partial \Omega^+}
u\partial_n v - v\partial_n u = 0.$$
The boundary conditions lead to:
\begin{equation}
  \label{eq:2}
  \int_{\Omega^+} \partial_t uv -\Delta v u - \int_{x=1}
u\partial_x v  + \int_{x=0} \partial_x u (0^+) v-u(x=0^+)\partial_x v =0.  
\end{equation}
We sum (\ref{eq:1}) and (\ref{eq:2}) :
\begin{multline*}
\int_{\Omega} \partial_t uv -\Delta v u - \int_{x=1}
u\partial_x v   + \int_{x=-1}
u\partial_n v - v(1-\alpha)\\+ \int_{x=0} v\left(\partial_x
  u(0^+)- \partial_x u (0^-)\right) +
\left(u(x=0^-)-u(x=0^+)\right)\partial_x v = (\alpha - \beta) \int_{\Omega^-} v.
\end{multline*}
We now make use of  the transmission conditions (\ref{eq:htransmissionb}) at $x=0$, leading to
\begin{multline*}
\int_{\Omega} \partial uv -\Delta v u - \int_{x=1}
u\partial_x v   + \int_{x=-1}
u\partial_n v - v(1-\alpha)\\- \int_{x=0} \beta v -\int_{x=0}
\alpha u(x=0^+)\partial_x v = (\alpha - \beta) \int_{\Omega^-} v.
\end{multline*}

Integrating by parts again, we finally obtain:

\begin{multline}
\label{eq:var+}
\int_{\Omega} \partial_t uv +\nabla v \nabla u -\int_{x=-1} v(1-\alpha)\\- \int_{x=0} \beta v -\int_{x=0}
\alpha u(x=0^+)\partial_x v - (\alpha - \beta) \int_{\Omega^-} v = 0,
\end{multline}
which is the weak formulation of our problem.

The corresponding equation is:
\begin{equation}
\label{eq:hmodel.3}
\begin{cases}
  -\Delta u + \partial_t u =  (\alpha-\beta){\mathbf 1}_{x<0} & \\
  -\alpha \partial_x\left(u(x=0^+)\delta_{x=0} \right) + \beta \delta_{x=0}&
  \text{in } \Omega = \{ -1<x<1\}, \\
  \partial_n u = 1- \alpha & \text{on} \{x=-1\}, \\
  \partial_n u = 0 & \text{on } \{x=-1\}, \\
  u \text{ is }1-\text{periodic in } y.
\end{cases}
\end{equation}

\begin{remarque}
  Taking $v=1$ in  (\ref{eq:var+}) leads to the energy conservation equation:
  $$\frac{d}{dt} \int_\Omega u = 1.$$
\end{remarque}

We compute the solution of (\ref{eq:var+}) using the numerical approximation of the Dirac mass at $x=0$:

$$  \int_{x=0} \alpha u(x=0^+)\partial_x v \approx \int_{\Omega}
\frac{\alpha}{\delta} \11_{0<x<\delta} u \partial_x v, $$
for $\delta$ arbitrary small. Note that this approximation may not be satisfactory as $\alpha$ increases. Taking a piecewise linear approximation of the Dirac mass around $x=0$ may improve the numerical treatment. This approach nevertheless diverges from the direct calculation for bigger $\alpha$ (see figure~\ref{fig:plot_zp6}), especially in the vicinity of the crack interface at $x=0$.

\subsection{Fixed-point approach}

Another numerical approach to solve (\ref{eq:hmodel.1})-(\ref{eq:hmodel.2}) that might be more accurate is to use the fixed-point approach, in the spirit of section 5 and remark 5.2. We apply the following iterative process:
\begin{itemize}
\item{Starting with an initial guess $F^0$ corresponding the flux imposed on $\Gamma^\beta$, we solve the problem (\ref{eq:homog_-}), and get the unique solution  $u^{+,0}$.}
\item{Then, we compute the trace $g^0$ of $(1-\alpha)u^{+,0}$ on $\Gamma^\beta$, and solve (\ref{eq:homog_+}) with the data $g=g^0$. This gives us 
the unique solution $u^{-,0}$.}
\item{We compute a new flux on $\Gamma^\beta$ $F^1 = \beta-\partial_n u^{-,0}$ and we repeat the first step with $F^1$.}
\end{itemize}

This procedure builds a converging sequence $(F^n)_{n\in\NN}$ in $L^2(\Gamma^\beta)$, together with the corresponding
solutions $u^{\pm,n}$ and the data $g^n$. We assess the convergence of $(F^n)_{n\in\NN}$ when the relative change between two successive iterations $\|\displaystyle\frac{F^{n+1}-F^{n}}{F^{n}}\|$ is smaller than a small fixed parameter. This has to be carried out at each time step of the simulation. 

\subsection{Synthesis}

\begin{figure}
        \centering
        \begin{subfigure}[b]{0.5\textwidth}
                \centering
                \includegraphics[width=\textwidth]{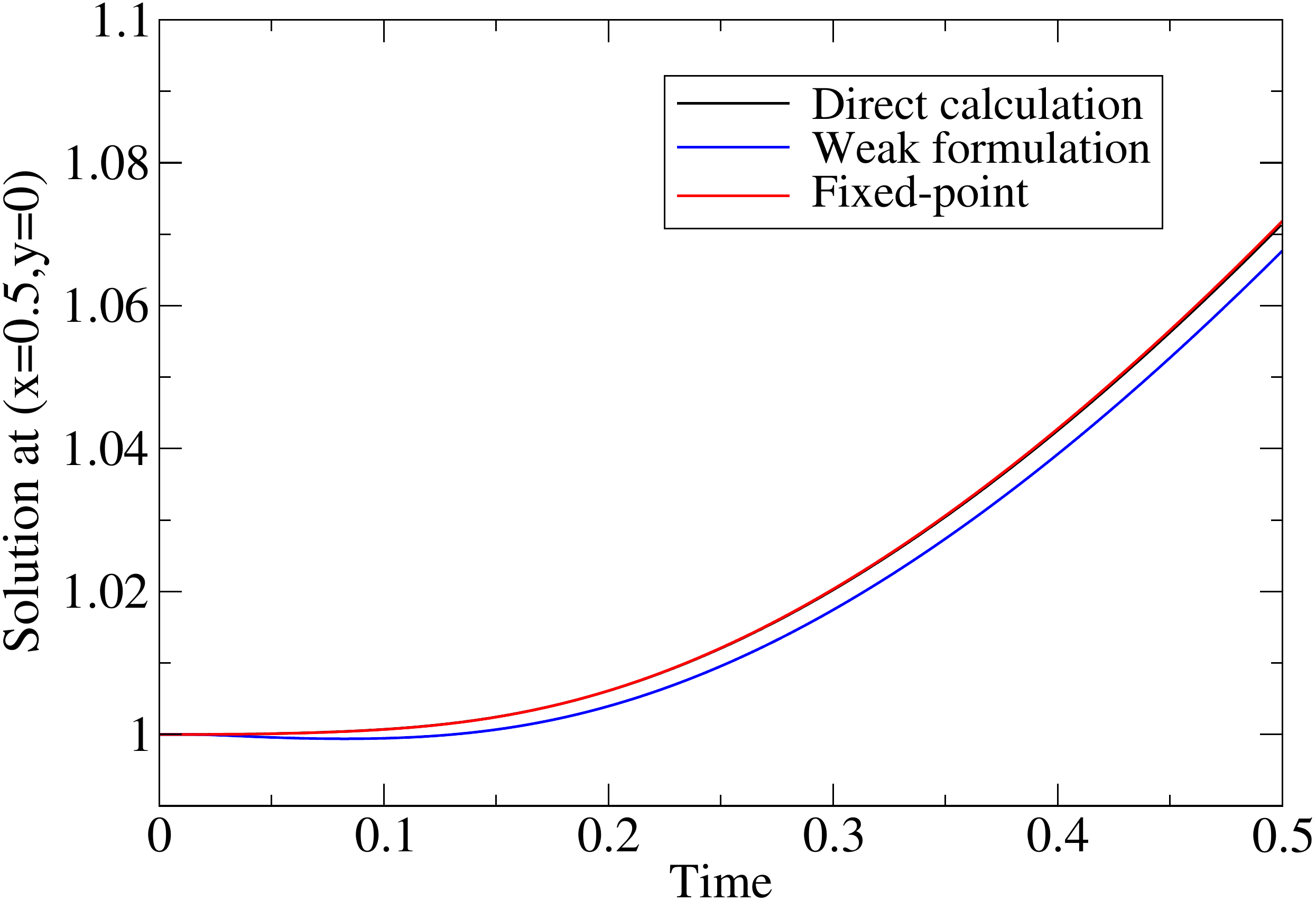}
                \caption{\textbf{Time evolution of $u(0.5,0)$}.}
                \label{fig:gull}
        \end{subfigure}%
        ~ 
        \begin{subfigure}[b]{0.5\textwidth}
                \centering
                \includegraphics[width=\textwidth]{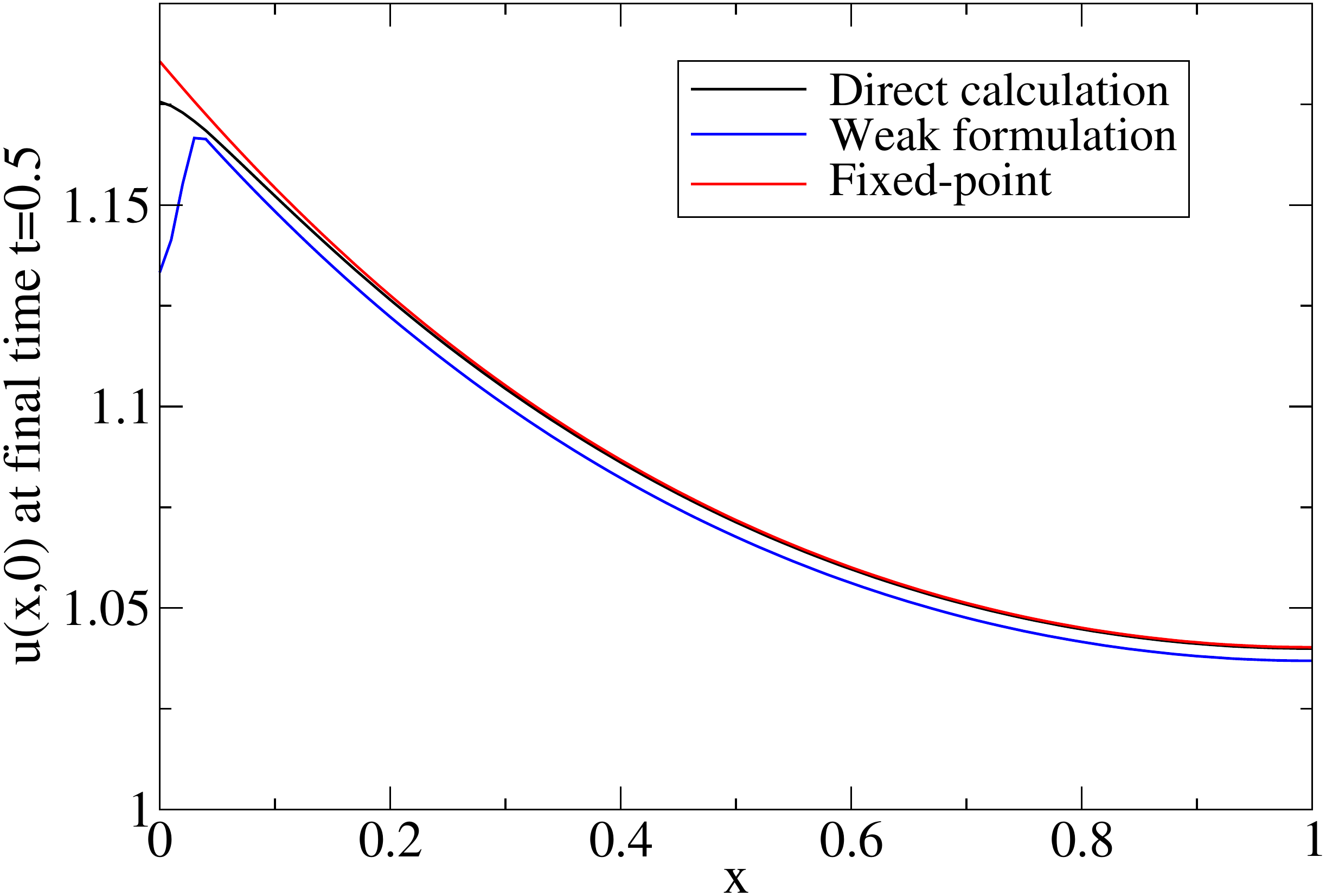}
                \caption{\textbf{Spatial profile at time $t_1=0.5$}.}
                \label{fig:tiger}
        \end{subfigure}
        \caption{\textbf{Comparison of the direct and homogenized approaches for $\alpha=0.1$}.\textit{The black curve is the direct calculation of the crack, i.e the limit solution of (\ref{eq:simple}) as $\varepsilon\to 0$. The two other curves correspond to the solution of the homogenized problem (\ref{eq:hmodel.1})-(\ref{eq:hmodel.2}) computed by two approaches: in red, it is the fixed-point method whereas the blue curve represents the numerical solution of the weak formulation (\ref{eq:var+}).}}\label{fig:plot_zp1}
\end{figure}

\begin{figure}
        \centering
        \begin{subfigure}[b]{0.5\textwidth}
                \centering
                \includegraphics[width=\textwidth]{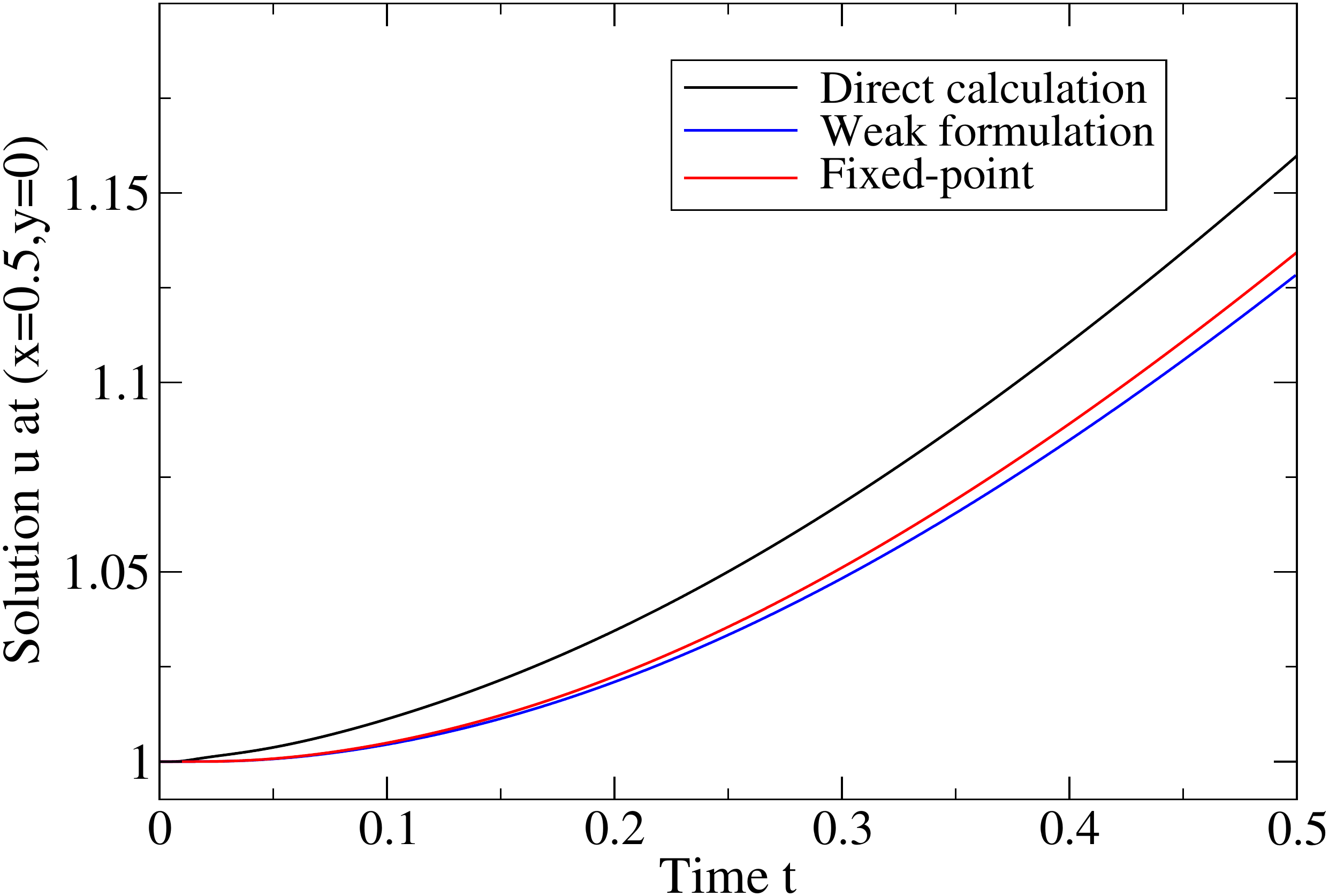}
                \caption{\textbf{Time evolution of $u(0.5,0)$}}
                \label{fig:gull2}
        \end{subfigure}%
        ~ 
        \begin{subfigure}[b]{0.5\textwidth}
                \centering
                \includegraphics[width=\textwidth]{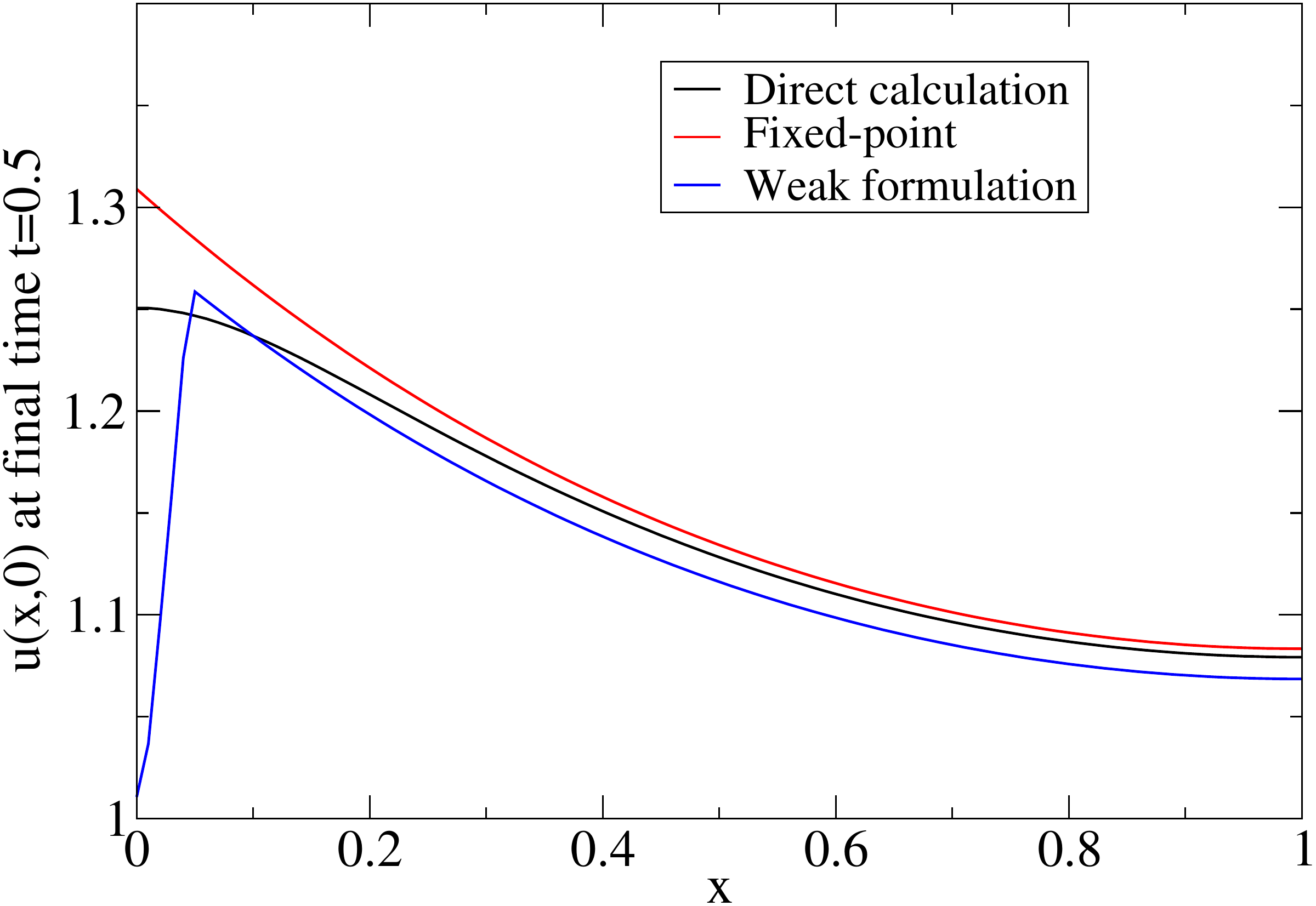}
                \caption{\textbf{Spatial profile at time $t_1=0.5$}}
                \label{fig:tiger2}
        \end{subfigure}
        \caption{\textbf{Comparison of the direct and homogenized approaches for $\alpha=0.6$}.\textit{We note that when we get close to $x=0$, i.e the bottom of the crack, the agreement between the different approaches is not so good. This is due to the singularity represented by a Dirac mass at $x=0$ in the weak formulation of the homogenized problem.}}\label{fig:plot_zp6}
\end{figure}

We thus compute the results given by each method:
\begin{itemize}
\item{direct calculation on the cracked domain $\Omega_{\varepsilon}$,}
\item{weak formulation on the homogenized domain $\Omega$,}
\item{fixed point approach on $\Omega$,}
\end{itemize}
for two crack configurations: $\alpha=0.1$ (small crack) and $\alpha=0.6$ (big crack). We fix $\beta=0$ in the simulations for the sake of simplicity.  We compare the different simulations by plotting the time evolution of the solution at $x=0.5$ (that is to say, in the core of the intact part of the cracked material) and the spatial profile of the solution at final time $t_1=0.5$. Results are shown in figures \ref{fig:plot_zp1} and \ref{fig:plot_zp6}. We can note that the fixed-point method appears to be more accurate than the solution of the weak formulation (\ref{eq:var+}) especially as $\alpha$ increases. This is due to the fact that (\ref{eq:var+}) involves a Dirac mass at $x=0$, proportional to the width of the crack $\alpha$. This term is treated approximately in our finite element simulations and leads to more significant errors for greater $\alpha$. In some way, the fixed-point method amounts to treat the Dirac mass at $x=0$ exactly.

\begin{figure}[h!]
  \centering
  \includegraphics[width=0.7\linewidth]{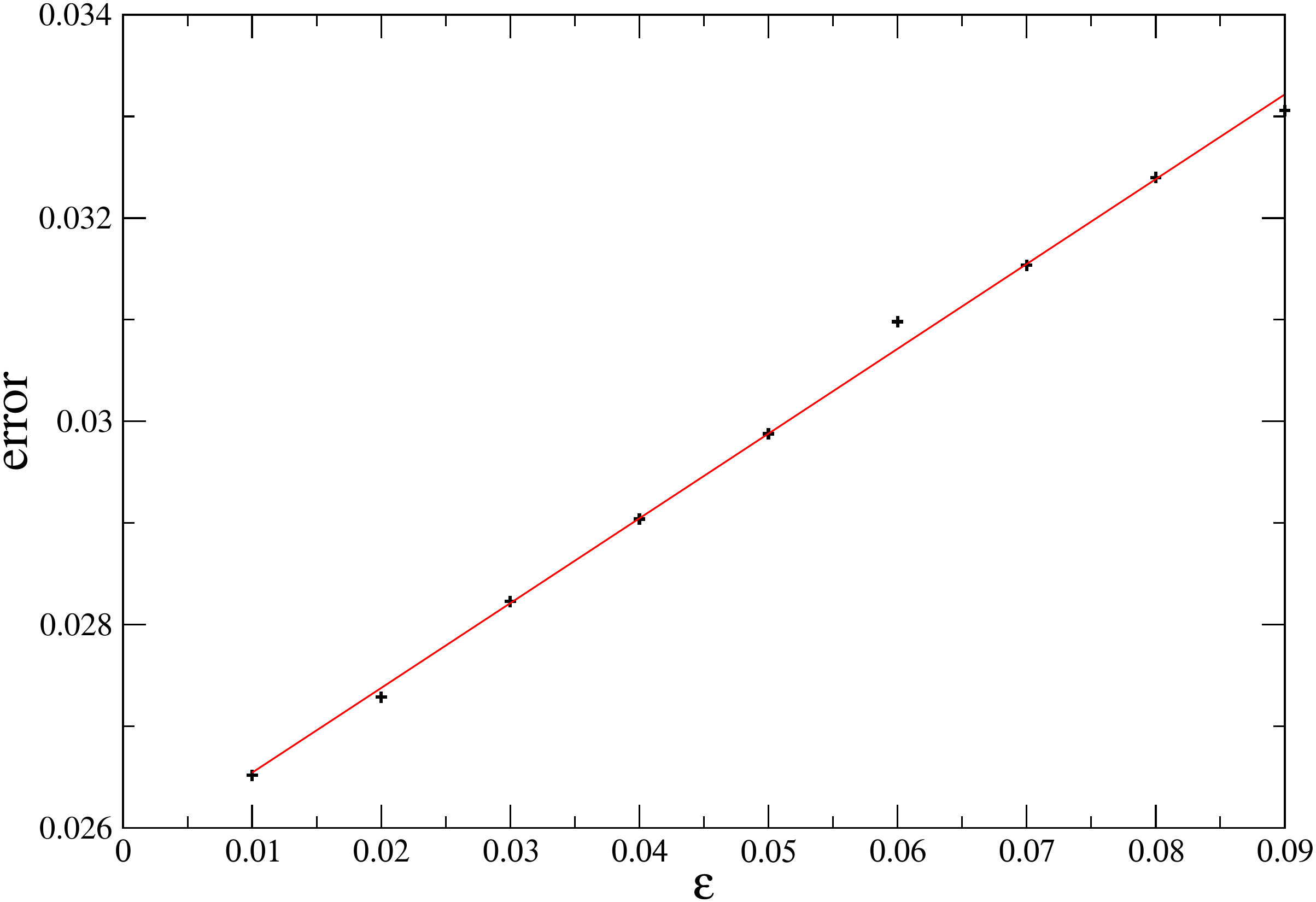}
  \caption{The error associated with the homogenized model as a function of $\varepsilon$.\textit{The error depends linearly on the period $\varepsilon$. As $\varepsilon\to 0$, the error tends to the residual error linked to the mesh size used in the finite element calculations, $\sim$ 0.02 here.} \label{fig:erreps}}
\end{figure}

We conclude this numerical illustration by studying the error associated with the homogenized model (\ref{eq:hmodel.1})-(\ref{eq:hmodel.2}) with respect to the period $\varepsilon$ of the cracked domain. More precisely, we calculate the error defined by:

\begin{equation}
  \label{eq:errl2}
\mbox{err} = \displaystyle\frac{\|\tilde{u}_\varepsilon-u_{\varepsilon}\|_{L^2(\Omega_\varepsilon^+)}}{\|u_{\varepsilon}\|_{L^2(\Omega_\varepsilon^+)}},
\end{equation}
where $u_\varepsilon$ is the solution of the exact model (\ref{eq:simple}) in the cracked domain $\Omega_\varepsilon$ and $\tilde{u}_\varepsilon$ is the projection of the solution $u$ of the homogenized model (\ref{eq:hmodel.1})-(\ref{eq:hmodel.2}) on the sub-domain $\Omega_\varepsilon^+$.
The error is plotted as a fonction of $\varepsilon$ on figure (\ref{fig:erreps}). It shows that the error depends linearly on the period $\varepsilon$. Besides, as $\varepsilon\to 0$, the error tends to the residual error linked to the mesh size used in the finite elements calculations.

\section{Conclusion and future work}

We have shown that a diffusion process initiated by an incoming flux through a periodic cracked medium can be modeled by a volume source term in the diffusion equation solved in the homogenized domain. The crack induces also a singularity at $x=0$ giving rise to some complications in the formulation of the homogenized problem. We have introduced a boundary layer around $x=0$ to treat it properly, and it leads to a Dirac mass located at the bottom of the crack in the weak formulation of the homogenized problem. Note that this singularity might be interpreted physically by considering that the temperature field in the vicinity of the crack is well described by a ponctual source located at the bottom. This remains true as long as we observe the temperature sufficently far from the crack, where we are not too sensitive to the details of the fracture profile. This can be shown rigorously in a very particular configuration for which the shape of the fracture is smooth enough to be described analytically using conformal mappings \cite{bui1,bui2}. 

Besides, we have developped our method in a very particular setting, where the cracks are supposed to be orthogonal to the surface of the material and periodically arranged. This is an elementary case on which we have tested our approach. We may now have in view to address much more general situations on which MOSAIC could be applied. For example, we may add some \textbf{stochastic} features, e.g the width of the crack $\alpha$ (which was supposed to be a fixed parameter) may become a stochastic variable. The periodicity of the structure may then disappear. The position of the crack might obey a stochastic process as well, and the crack setting may thus become very general. It would be interesting to see how the model developed here and proved on an elementary pattern may be extended to more intricate stochastic configurations. Even more generally, the homogenization of a diffusion process through fractal cracks, in the spirit of \cite{achdou}, may be a challenging question to tackle.

\begin{figure}[h!]
  \centering
  \includegraphics[width=0.7\linewidth]{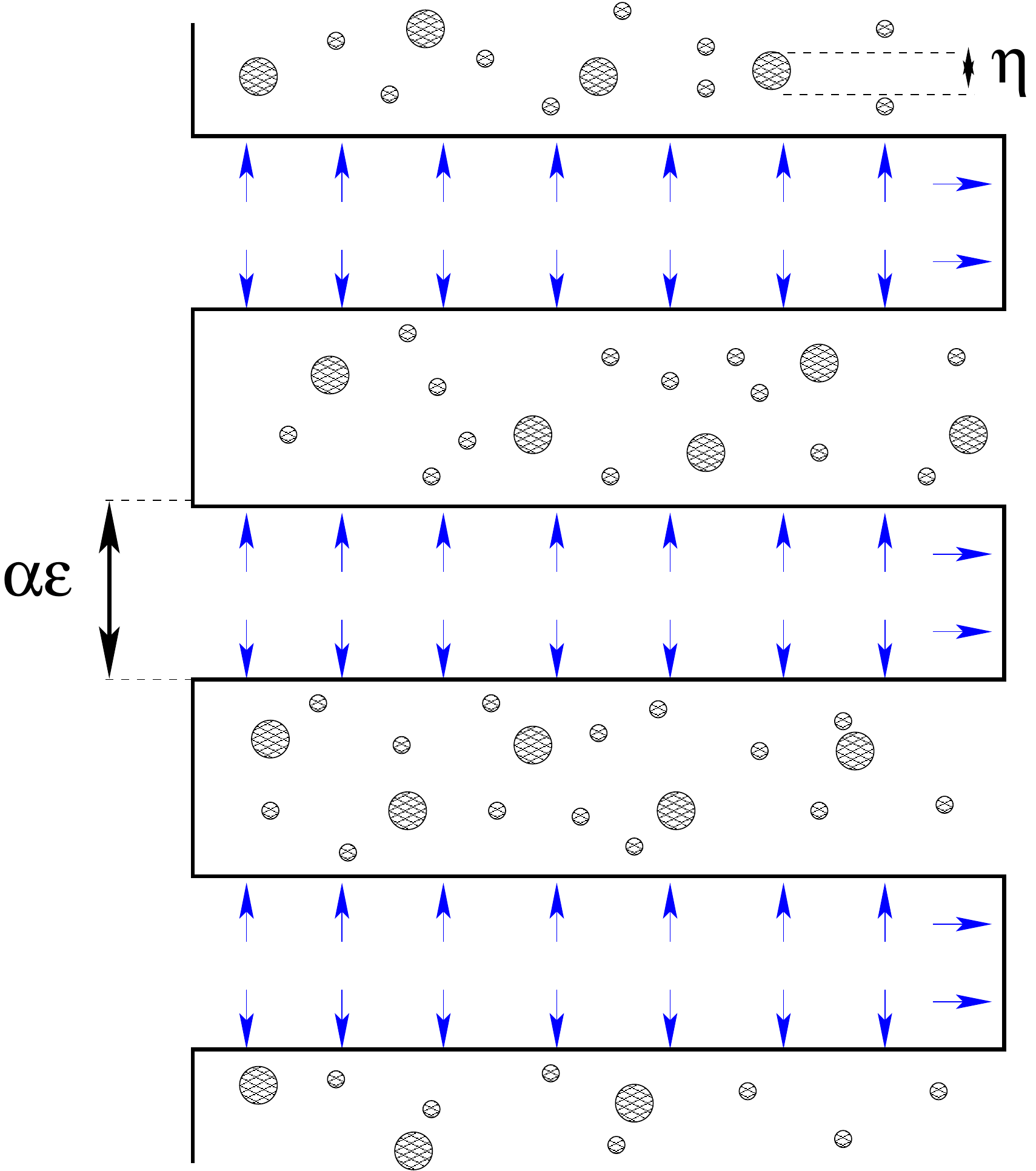}
  \caption{An example of a cracked heterogeneous medium.\textit{The medium is characterized by two length scales: $\alpha\varepsilon$ corresponding to the width of the fractures on which an incoming flux is imposed (symbolized by the blue arrows), and $\eta << \varepsilon$ which is the scale of smaller heterogeneities. The fractures of width $\alpha\varepsilon$ can be treated by the MOSAIC approach whereas the smaller structures may be described by an effective approach.} \label{fig:subcrack}}
\end{figure}

Another extension of the present work may be to consider that the cracked material is not homogeneous but filled with radiation-free micro-cracks (see figure \ref{fig:subcrack} for an illustration). Provided the size of such micro-cracks is much smaller than $\alpha$, the effective behaviour of the micro-cracked medium can be estimated - and in some cases rigorously bounded - using general methods from the theory of composites materials, possibly extended in a non-linear regime (\cite{tw85}, \cite{ppc91}, \cite{mp05}). Those methods make use of the so-called "translation method" devised by Murat and Tartar \cite{murat}. Combining such approaches with the MOSAIC method is the object of ongoing work  \cite{mb}.
 
\begin{thebibliography}{2}
\bibitem{achdou} Y. Achdou, C. Sabot, N. Tchou, \textbf{Diffusion and propagation problems in some ramified domains with a fractal boundary} , M2AN Math. Model. Numer. Anal. 40 (2006), no. 4, 623--652




\bibitem{blp} A. Bensoussan, J.L Lions, G. Papanicolaou, \textbf{Asymptotic analysis of periodic structures}, Studies in Mathematics and its Applications, 5. North-Hollan Publishing Co., Amsterdam-New York, 1978.

\bibitem{dautray-lions} Dautray, Robert and Lions, Jacques-Louis and Artola, Michel and Bardos,
    Claude and Cessenat, Michel and Kavenoky, Alain and Lascaux, Patrick and
    Mercier, Bertrand and Pironneau, Olivier and Sentis, R\'emi,
{\bf Analyse math\'ematique et calcul num\'erique pour les sciences et les
    techniques. Volume 9: \'Evolution: num\'erique, transport. (Mathematical
    analysis and numerical methods for science and technology. Volume 9:
    Evolution: numerical methods, transport equations).} 
    , INSTN Collection Enseignement. Paris etc.: Masson. xxxiv, 1988.

  \bibitem{ding} Z. Ding, \textbf{A proof of the trace theorem of Sobolev spaces on Lipschitz domains.}
Proc. Amer. Math. Soc. 124 (1996), no. 2, 591--600.

\bibitem{evans} L. C. Evans, \textbf{Partial Differential equations,}
  Graduate Studies in Mathematics, 19. American Mathematical Socieity,
  Providence, RI, 2010.

\bibitem{gagliardo} E. Gagliardo, \textbf{Caratterizzazioni delle tracce sulla frontiera relative ad alcune classi di funzioni in n variabili.} (Italian)
Rend. Sem. Mat. Univ. Padova 27 1957 284--305. 

\bibitem{free07} F. Hecht, O. Pironneau, A. Le Hyaric, K Ohtsuke,
  FreeFem++ (manual). http://www.freefem.org, 2007.

\bibitem{lions-magenes} J.-L. Lions, E. Magenes, \textbf{
Probl\`{e}mes aux limites non homog\`{e}nes et applications.} Vol. 1. (French)
Travaux et Recherches Math\'{e}matiques, No. 17 Dunod, Paris 1968.

\bibitem{mih99} Mihalas and Mihalas, \textbf{Foundations of Radiation Hydrodynamics}, Dover 1999.

\bibitem{bui1} Bui H.D., Erlacher A., \textbf{Propagation dynamique d'une zone endommag\'{e}e
dans un solide \'{e}lastique fragile en mode III et en r\'{e}gime permanent},
C.R. Acad Sci Paris, 290, p.273-276, 1980


\bibitem{bui2} Bui H.D., Erlacher A., Nguyen QS, \textbf{Propagation de fissure en thermo\'{e}lasticit\'{e}
dynamique}, J.M\'{e}canique, Vol 19, No.4, 1980.

\bibitem{tw85} Talbot D.R.S., Willis J.R.,  \textbf{Variational principles for inhomogeneous nonlinear media}, IMA J. Appl. Math, 35, 39-54,1985.

\bibitem{ppc91} Ponte Castaneda P., \textbf{The effective mechanical properties of nonlinear isotropic composites}, J. Mech. Phys. Solids. 39, 45-71,1991.

\bibitem{mp05} Peigney M., \textbf{A pattern based method for bounding the effective response of a nonlinear composite}, J. Mech. Phys. Solids, 53, p.923-948, 2005.

\bibitem{murat} Murat F., Tartar L., \textbf{Calcul des variations et homog\'{e}n\'{e}isation} in \textit{Les M\'{e}thodes de l'Homog\'{e}n\'{e}isation: Th\'{e}orie et Applications en Physique}, cours de l'Ecole d'\'Et\'{e} d'Analyse Num\'{e}rique ,Eyrolles, Paris,1985.

\bibitem{mb} Peigney M., Peigney B., in preparation.
\end{thebibliography}
\end{document}